\documentclass[12pt]{article}
\pdfoutput=1

\usepackage{url}
\makeatletter
\newcommand{\address}[1]{\gdef\@address{#1}}
\newcommand{\email}[1]{\gdef\@email{\url{#1}}}
\newcommand{\sites}[1]{\gdef\@sites{\url{#1}}}
\newcommand{\@endstuff}{\par\vspace{\baselineskip}\noindent\small
	\begin{tabular}{@{}l}\scshape\@address\\\textit{E-mail address:} \@email \\ \textit{URL:} \@sites \end{tabular}}
\AtEndDocument{\@endstuff}
\makeatother

\usepackage{titling}
\title{\normalsize\textbf{{\large I}NTERMEDIATE {\large G}EODESIC {\large G}ROWTH IN \\ {\large V}IRTUALLY {\large N}ILPOTENT {\large G}ROUPS} \vspace*{-5mm}}
\author{Corentin Bodart \vspace*{-5mm}}
\date{\today}

\address{Mathematical Institute, University of Oxford, UK}
\email{corentin.bodart@maths.ox.ac.uk}
\sites{https://sites.google.com/view/corentin-bodart}

\usepackage{import}
\usepackage[backend=bibtex, style=alphabetic, sorting=nyt,
    doi=false,isbn=false,url=false,eprint=false, maxnames=50]{biblatex}
\usepackage[english]{babel}
\usepackage{amsmath}
\usepackage{amsfonts}
\usepackage{amssymb}
\usepackage{mathrsfs} 

\DeclareMathAlphabet{\funcal}{U}{BOONDOX-cal}{m}{n}
\SetMathAlphabet{\funcal}{bold}{U}{BOONDOX-cal}{b}{n}
\DeclareMathAlphabet{\funbcal} {U}{BOONDOX-cal}{b}{n}

\usepackage{dsfont}

\usepackage[dvipsnames]{xcolor}

\usepackage{hyperref}
\usepackage{enumitem}
\usepackage[chapter]{easy-todo}         

\usepackage[left=3cm,right=3cm,top=2.8cm,bottom=3cm, marginparsep=3mm]{geometry}
\usepackage{changepage}

\usepackage[retainorgcmds]{IEEEtrantools}   
\usepackage[font=small, hypcap=false]{caption}
\usepackage{parskip}

\usepackage{tikz-cd}
\usepackage{tikz}
\usetikzlibrary{shapes, shapes.geometric, arrows, arrows.meta, decorations.markings, decorations.pathreplacing}
\pgfdeclarelayer{background} 
\pgfdeclarelayer{foreground}
\pgfsetlayers{background,main,foreground}
\import{Packages/}{custom_macros.tex}
\usepackage{xcolor}

\usepackage{amsthm}
\usepackage{chngcntr}

\theoremstyle{plain}
\newtheorem{thm}[equation]{Theorem}
\newtheorem*{thm*}{Theorem}
\newtheorem{lemma}[equation]{Lemma}
\newtheorem*{lemma*}{Lemma}
\newtheorem{cor}[equation]{Corollary}
\newtheorem*{cor*}{Corollary}
\newtheorem{prop}[equation]{Proposition}
\newtheorem*{prop*}{Proposition}

\theoremstyle{definition}
\newtheorem{defi}[equation]{Definition}
\newtheorem{rem}[equation]{Remark}
\newtheorem*{rem*}{Remark}

\begin{document}

\maketitle

\begin{abstract}
   We give a criterion on pairs $(G,S)$ - where $G$ is a virtually $s$-step nilpotent group and $S$ is a finite generating set - saying whether the geodesic growth is exponential or strictly sub-exponential. Whenever $s=1,2$, this goes further and we prove the geodesic growth is either exponential or polynomial. For $s\ge 3$ however, intermediate growth is possible. We exhibit a pair $(G,S)$ for which $\gamma_\geod(n) \asymp \exp\!\big(n^{3/5}\cdot \log(n)\big)$, where $G$ contains a $3$-step nilpotent group - the Engel group - as a finite-index subgroup. This is the first example of group with intermediate geodesic growth. Along the way, we prove results on the geometry of virtually nilpotent groups, including asymptotics with error terms for their volume growth, and disprove a conjecture by Breuillard and Le Donne.
   \medbreak
    \vspace*{-2mm}
\end{abstract}

A common scheme in Geometric Group Theory is to consider groups as geometric spaces, define some notion of growth and then classify groups with growth in a given regime by their algebraic properties. The first example is volume growth: we consider a group $G$ together with a finite (monoid) generating set $S$. This defines a word metric
\[ \norm g_S =  \min \{ \ell(w) \mid w\in S^* \text{ and }\overline w = g \}\]
on $G$. The \emph{volume growth function} of the pair $(G,S)$ is
\[ \beta_{(G,S)}(n) = \#\!\left\{ g\in G \;\big|\; \norm{g}_S\le n \right\}. \]
Milnor famously asked two questions about the volume growth of groups \cite{Milnor}:
\begin{itemize}[leftmargin=6mm, label=\textbullet]
	\item For which groups do we have $\beta_{(G,S)}(n)\preceq n^d$ for some constant $d\ge 0$? Milnor further suggested the family of virtually nilpotent groups.
	\item Does the volume growth necessarily satisfy $\beta_{(G,S)}(n)\preceq n^d$ or $\beta_{(G,S)}(n)\succeq \exp(n)$?
\end{itemize}
On the one hand, Gromov confirmed that groups with polynomial volume growth coincide with virtually nilpotent groups. Doing so, he introduced a large machinery to make the link with nilpotent Lie groups \cite{Gromov}. On the other hand, Grigorchuk constructed a family of groups of intermediate growth (i.e., whose growth is neither bounded above by polynomials, nor below by exponentials). This family is now source of groups with many intriguing properties \cite{Grigorchuk_growth}. These results constitute an ideal to aim for.

\newpage

In our case, the story starts in the same way: a pair $(G,S)$, and the word metric $\norm{\,\cdot\,}_S$. However, instead of counting elements, we count geodesics. A word $w\in S^*$ is a \emph{geodesic} if $\ell(w)=\norm{\overline w}_S$. The \emph{geodesic growth function} of $(G,S)$ is then defined as
\[ \gamma(n) = \gamma_\geod^{(G,S)}(n) = \#\!\left\{ w\in S^* \;\big|\; w\text{ is geodesic and }\ell(w)\le n \right\}.\]

We should mention that geodesic growth is sensitive to the choice of a generating set. Indeed, Bridson, Burillo, Elder and {\v{S}}uni{\'c} proved that any infinite group $G$ admits a generating set $S$ such that $\gamma_\geod^{(G,S)}(n)$ grows exponentially \cite[Example 6]{bridson2012groups}.\footnote{This holds if we allow generating \emph{multisets}. Otherwise, the only virtually nilpotent counter-example is $G=\ZZ$ (as a corollary of Theorem \ref{thm:crit}), and any other hypothetical counter-example would be of intermediate volume growth (as exponential volume growth implies exponential geodesic growth). \vspace*{1mm}} Therefore, Bridson \emph{et al}.\ ask the following questions:
\begin{adjustwidth}{3mm}{3mm}
	\textbf{Question 1.} Characterize groups $G$ with polynomial geodesic growth, that is, with $\gamma_\geod^{(G,S)}(n)\preceq n^d$ for some constant $d\ge 0$ and for at least one generating set $S$.
	
	\textbf{Question 2.} Does there exists a pair $(G,S)$ with intermediate geodesic growth? \vspace*{1mm}
\end{adjustwidth}

In the same paper, Bridson \emph{et al}.\ proposed some partial answers. Their main theorem is a sufficient condition for polynomial geodesic growth:
\begin{thm*}[{\cite[Theorem 1]{bridson2012groups}}]
	Let $G$ be a finitely generated group. If there exists an element $a\in G$ such that $H=\lla a\rra_G$ is a finite-index abelian subgroup, then there exists a symmetric generating set $S$ such that $(G,S)$ has polynomial geodesic growth.
\end{thm*}

This includes all virtually cyclic groups, and also groups like
\[ \funcal{v\!Z} = \ZZ^2 \rtimes C_2 = \la a,t \mid t^2=[a,a^t]=1 \ra, \]
where $C_2=\la t\mid t^2=e\ra$. Subsequently, Bishop and Elder \cite{bishop2020virtHeis} proved that the group
\[ \funcal{v\!H} = H_3(\ZZ)\rtimes C_2 = \la a,t \mid t^2=[a,[a,a^t]]=1 \ra \]
has polynomial geodesic growth w.r.t.\ the generating set $S=\{a^\pm,t\}$.

In the opposite direction, Bridson \emph{et al}.\ showed that any group factoring onto $\ZZ^2$ has exponential geodesic growth w.r.t.\ every generating set. This applies for finitely generated nilpotent groups which are not virtually cyclic \cite[Lemma 13]{bridson2012groups}. We generalize most of these results in the following criterion.

\textbf{Setup.}
	Let $G$ be a virtually $s$-step-nilpotent group, with $S$ a finite generating set. Consider $H$ a torsion-free, $s$-step nilpotent, finite-index, normal subgroup of $G$.\footnote{We can always find such an $H$: $G$ has a finite-index $s$-step nilpotent group $H'$, which contains a finite-index torsion-free subgroup $H''$ \cite[Thm 3.23]{torsionfree}. Finally $H''$ contains a finite-index subgroup
	\[ H = \mathrm{core}(H'') = \bigcap_{g\in G}gH''g^{-1}\len G.\]} \vspace*{2mm}
	
	\newpage
	
	We consider the \emph{isolator} of the commutator subgroup
	\[ [[H,H]] = \big\{ h\in H \mid \exists n>0,\, h^n\in[H,H]\big\}. \]
	We get out of this data a map $\pi\colon H\onto H/[[H,H]]\simeq\ZZ^d$, and an action of the finite group $F=G/H$ on $H/[[H,H]]$ (by conjugation). We define the multiset\footnote{$A(S)$ is multiset in the following sense: any point $p\in\QQ^d$ appears as many times in $A(S)$ as there are simple cycles $a\in S^*$ such that $p=\pi(\bar a)/\ell(a)$.}
	\[ A=A(S) = \left\{ \frac{\pi(\bar a)}{\ell(a)} \in \QQ^d \;\middle|\; a\in S^*\text{ labels a simple cycle in }\Sch( H\backslash G, S)\right\}, \]
	where $\Sch( H\backslash G, S)$ denotes the Schreier graph. Finally, we define a polytope $P(S)=\Conv(A(S)^F)$, where $A^F$ denotes the orbit of $A$ under conjugation by $F$.
	
\begin{thm} \label{thm:crit} If no two elements of $A(S)$ lie on a common facet of $P(S)$, then the geodesic growth is subexponential. More precisely, we give the following upper bounds:
	\begin{itemize}[leftmargin=8mm, label={\normalfont\textbullet}]
		\item If $s\le 2$, the geodesic growth is bounded above by a polynomial.
		\item If $s\ge 3$, the geodesic growth is bounded above by
		\[ \gamma_\geod^{(G,S)}(n) \preceq \exp\!\big(n^{\alpha_s}\log(n)\big), \]
		with $0<\alpha_s<1$ an explicit constant {\normalfont(}eg. $\alpha_3=3/5${\normalfont)}.
	\end{itemize}
	Otherwise the geodesic growth is exponential.
\end{thm}

\bigbreak

In particular, we reduce the \say{virtually $2$-step nilpotent} case of Question 1 to the characterization of finite subgroups $F\le\GL_d(\ZZ)$ with a certain property, as follows:
\begin{cor} \label{thm:crit_step2}
Let $G$ be a finitely generated, virtually $2$-step-nilpotent group. Consider $H\len G$ a torsion-free, $2$-step nilpotent, finite-index, normal subgroup of $G$. This defines an action of $F=G/H$ on $\ZZ^d\simeq H/[[H,H]]$. The following assertions are equivalent:
\begin{enumerate}[leftmargin=8mm, label={\normalfont(\roman*)}]
	\item There exists a generating set $S$ such that $(G, S)$ has polynomial geodesic growth.
	\item There exists a finite set $A\subset \ZZ^d$ such that $P=\Conv(A^F)$ is a full-dimensional polytope, and no two elements of $A$ lie on the same facet of $P$.
\end{enumerate}
\end{cor}
\begin{proof}
	We have to prove (ii) $\Rightarrow$ (i). We construct a generating set from $A=\{p_1,\ldots,p_m\}$ as follows. Consider $S_0$ a fixed generating set for $G$, and define
	\[ S_n = S_0 \cup \{a_1^n,\ldots,a_m^n\} \text{ for all }n\in\ZZ_{>0}, \]
	where $a_i$ are elements of $H$ satisfying $\pi(a_i)=p_i$. Observe that the new generators $a_i^n\in H$ only add loops in the Cayley graph of $G/H$, and therefore
	\[A(S_n) = A(S_0)\cup nA. \]
	For $n$ large enough, we have $P(S_n)=\Conv(nA^F)$ and $A(S_0)\cap \partial P(S_n)=\emptyset$. At this point, the hypothesis on $A$ implies that $(G,S_n)$ has polynomial geodesic growth.
\end{proof} \vspace*{0mm}

\newpage

\begin{rem}
	The statement of this last result can be adapted if we only allow \emph{symmetric} generating sets $S$. The only other modification needed is
	\begin{itemize}[leftmargin=8mm]
		\item[\hypertarget{symm_cond}{(ii')}] There exists a \emph{symmetric} finite set $A\subset \ZZ^d$ such that $P=\Conv(A^F)$ is a full-dimensional polytope, and no two elements of $A$ lie on the same facet of $P$.
	\end{itemize}
	Note that conditions (ii) and (ii') are not equivalent. An example is given by the group $G_2=\ZZ^2\rtimes C_2$  where $C_2=\la r\ra$ acts by $180^\circ$ rotations (see also \cite[Example 16]{bridson2012groups}). If we only look at \emph{symmetric} sets $A$, we always have $A^{C_2}=A\cup -A=A$, so that both vertices of any facet of $P$ belong to $A$. In contrast $G_2$ satisfies condition (ii):
	\begin{center}
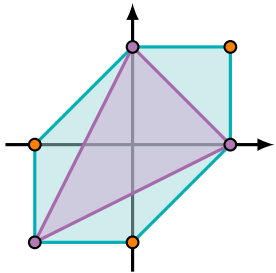

		\begin{tikzpicture}[scale=1.3]
			\draw[very thick, -latex] (-1.3,0) -- (1.45,0);
			\draw[very thick, -latex] (0,-1.3) -- (0,1.45);
			\draw[very thick, TealBlue, fill=TealBlue!40, fill opacity=.5] (1,0) -- (1,1) -- (0,1) -- (-1,0) -- (-1,-1) -- (0,-1) -- cycle;
			\draw[thick, black, fill=TealBlue] (1,1) circle (1.3pt);
			\draw[thick, black, fill=TealBlue] (-1,0) circle (1.3pt);
			\draw[thick, black, fill=TealBlue] (0,-1) circle (1.3pt);
			\draw[thick, black, fill=Purple!70] (1,0) circle (1.7pt);
			\draw[thick, black, fill=Purple!70] (0,1) circle (1.7pt);
			\draw[thick, black, fill=Purple!70] (-1,-1) circle (1.7pt);
		\end{tikzpicture} \vspace*{-3mm}
		\captionof{figure}{$A=\{(1,0),(0,1),(-1,-1)\}$ in purple and $P(S)$ in green.}
	\end{center}
	 This means that the geodesic growth of $G_2$ is polynomial w.r.t.\ $S=\{x,y,(xy)^{-1},r\}$, and exponential w.r.t.\ any symmetric generating set (as shown in \cite{bridson2012groups}).
\end{rem}

\bigbreak

Regarding Question 2, we get the following affirmative answer:
\begin{thm} \label{thm:inter}
	The geodesic growth of the group
	\[ \vEc = \la a,t \;\big|\; t^2=1;\; [a,[a,a^t]]=[a^t,[a,a^t]] \;\text{\normalfont commutes with }a,a^t\ra \]
	with generating set $S=\{a^{\pm 1},t\}$ satisfies $\gamma_\geod(n) \asymp \exp\!\big(n^{3/5}\cdot\log(n)\big)$.
\end{thm}

Note that this group is virtually $3$-step nilpotent, more precisely it admits an index-$2$ subgroup isomorphic to the so-called \say{Engel group}. In some sense, this is the next smallest candidate for intermediate geodesic growth, as virtually abelian have either polynomial or exponential geodesic growth (see \cite{bishop2021virtAbel}), and the same holds true in virtually $2$-step nilpotent groups (see Theorem \ref{thm:crit}). The construction relies on the same trick as examples $\funcal{v\!Z}$ and $\funcal{v\!H}$ of Bridson \emph{and al}.\ and Bishop-Elder. Our proof re-uses some of their ideas, combined with insights from nilpotent geometry.

\bigbreak

As a byproduct, we provide estimates on the volume growth $\beta_{(G,S)}(n)$ of virtually nilpotent groups $G$. A celebrated result due to Pansu \cite{Pansu} states that
\[ \beta_{(G,S)}(n) = c_{(G,S)}\cdot n^d + o(n^d)\] 
whenever $G$ is virtually nilpotent. Subsequently the error term was refined whenever $G$ is $2$-step nilpotent \cite{Stoll_metric} and more generally $s$-step nilpotent \cite{BreuillardLeDonne, Gianella}. We extend these results to virtually $s$-step nilpotent groups:

\begin{thm}[Corollary \ref{cor:vol}]
	Let $G$ be a virtually $s$-step nilpotent group, and $S$ a finite symmetric generating set. The volume growth satisfies
	\[ \beta_{(G,S)}(n) = c_{(G,S)}\cdot n^d + O(n^{d-\min\{1-\alpha_s,\delta_s\}}),\]
	where $\alpha_s<1$ is defined in Section \ref{sec:control}, $\delta_s=1$ for $s=1,2$ and $\delta_s=\frac1s$ for $s\ge 3$.
\end{thm}
Finally, in Corollary \ref{rem:BrLD}, we disprove a conjecture by Breuillard and Le Donne stating that, in any torsion-free nilpotent group $H$ with a generating set $X$, we should have
\[ \norm{g}_X - \norm{g}_{\Stoll,X} = O_{H,X}(1) \]
(where $\norm{\,\cdot\,}_{\Stoll,X}$ is the Stoll metric, defined in Section \ref{sec:prelim}) \cite[Conjecture 6.5]{BreuillardLeDonne}. This conjecture was made in an effort to improve the volume growth estimate to $O(n^{d-1})$. 

\medbreak

\textbf{Organization of the paper.} In Section \ref{sec:prelim}, we recall some results mainly due to Stoll on what is now called the \say{Stoll metric}. Section \ref{sec:gen} gathers some results on word metrics in virtually nilpotent group. We compare word metrics in virtually nilpotent groups and the corresponding nilpotent subgroups, and give some local structure for geodesics. Section \ref{sec:main} is devoted to the proof of Theorem \ref{thm:crit}. In Section \ref{sec:Eng}, we do a deeper dive into the Engel group. We introduce a model for this group, provide fine lower bounds on word length of some elements, and prove the lower bound on geodesic growth needed for Theorem \ref{thm:inter}. Finally, Section \ref{sec:Qs} compiles some remarks and questions.

\medbreak

\textbf{Acknowledgements.} I'd like to thank Alex Bishop for helpful conversations, feedback, and suggesting a simplification for Proposition \ref{prop:back_from_depths}, Tatiana Nagnibeda for her constant support and feedback on early versions, and the referee for numerous suggestions. Special thanks go to Samuel Borodi for pointing out a mistake in the published version of the article. The author acknowledges support of the Swiss NSF grant 200020-200400.

\counterwithin{equation}{section}

\section{Preliminaries} \label{sec:prelim}

In this paper, we will consider metrics slightly more general than the usual word metrics on groups, allowing each generator to have a different weight:
\begin{defi}
	Consider $H$ a group, together with $X$ a finite set generating $H$ as a monoid ($X$ could be asymmetric, and could be a multiset) and a weight function $\sigma\colon X\to \ZZ_{>0}$. All words $w\in X^*$ can be written as $w=x_1^{m_1}x_2^{m_2}\ldots x_k^{m_k}$ for some letters $x_i\in X$ satisfying $x_i\ne x_{i+1}$, and $m_i\in\mathbb Z_{>0}$. We define
	\begin{itemize}[leftmargin=8mm]
		\item the \emph{coarse length} $k(w)=k$,
		\item the \emph{length} $\ell_\sigma(w) = \sum_{i=1}^k m_i\cdot \sigma(x_i)$.
	\end{itemize}
	Moreover, for each element $g\in G$, we define
	\[ \norm{g}_{X,\sigma} = \min\big\{ \ell_\sigma(v) \mid v\in X^* \text{ and } \bar v=g \big\}. \]
	Whenever $\sigma\equiv 1$, we will drop the subscript $\sigma$.
\end{defi}

\subsection{The Stoll metric}

Consider $H$ a torsion-free nilpotent group and $\bar H$ its Malcev completion. In the Malcev completion, exponentiation like $h^\mu$ makes sense for $h\in\bar H$ and $\mu\in\RR$.

\begin{defi}[$\RR$-words]
Fix a finite Lie generating set $X$ for $\bar H$.  An \emph{$\RR$-word} is an expression $w = x_1^{\mu_1}\cdot x_2^{\mu_2}\cdot\ldots\cdot x_k^{\mu_k}$ with $x_i\in X$ and $\mu_i\in\RR_{> 0}$. We denote the set of $\RR$-words by $X_\RR^*$. The notions of coarse length $k(w)$ and length $\ell_\sigma(w)$ extend to $X_\RR^*$.
\end{defi}
In \cite{Stoll_metric}, Stoll proves that any element $g\in\bar H$ can be represented by an $\RR$-word. We can therefore define the  Stoll norm on $\bar H$ as follows:
\begin{defi}[Stoll metric] Given $\bar H$ a simply connected nilpotent Lie group, $X$ a finite Lie generating set and $\sigma\colon X\to \ZZ_{>0}$ a weight function, we define \vspace*{-.5mm}
	\[ \norm{h}_{\Stoll,X,\sigma} = \inf\left\{ \ell_\sigma(w) \;\Big|\; w\in X_\RR^* \text{ and }\overline w= h \right\}. \vspace*{-2mm}\]
\end{defi}
Observe that $X^*\subset X^*_\RR$, therefore we have $\norm{h}_{X,\sigma}\ge \norm{h}_{\Stoll,X,\sigma}$ for all $h\in H$.
\subsection{The Abelian case}\label{sec:abelian}

We first have a look at the Abelian case and make some useful observations. See \cite{duchin2012geometry} for a more thorough treatment. We have $H\simeq \ZZ^d$ and $\bar H=\RR^d$. In this case, the Stoll metric has a very geometrical interpretation: 
\begin{lemma}\label{sec1:abelian} The Stoll metric coincides with a \say{Minkowski norm}: for all $v\in\RR^d$,
	\[ \norm{v}_{\Stoll,X,\sigma} = \norm{v}_{\Mink,P} \overset{\mathrm{def}}= \min\left\{\lambda\ge 0 \;\middle|\;v \in \lambda P\right\}, \]
	where $P = \Conv\big\{\frac{x}{\sigma(x)} \;\big|\; x \in X\big\}  \subset \RR^d$. Moreover, an $\RR$-word $x_1^{\mu_1}\cdot\ldots\cdot x_k^{\mu_k}$ is geodesic if and only if all $\frac{x_i}{\sigma(x_i)}$ lie on a common face of $P$.
\end{lemma}
\begin{proof} We can reduce ourselves to the case $\tilde\sigma\equiv 1$ by setting $\tilde X=\big\{\frac{x}{\sigma(x)}\mid x\in X\big\}\subset \bar H$. The case $v=0$ is trivial. Suppose $v\ne 0$ and let $m=\norm{v}_{\Mink,P} > 0$. 
	
	We first construct an $\RR$-word representing $v$ of length $m$. Consider $F$ the minimal face of $P$ containing $\frac1m \cdot v$. By the Caratheodory theorem, there exists $d$ vertices of $F$, say $x_1,\ldots,x_d\in X$, such that $\frac1m\cdot v \in \Conv(x_1,\ldots,x_d)$, i.e., \vspace*{-2mm}
	\[  \exists \nu_1,\ldots,\nu_d\ge 0 \quad\text{such that}\quad \nu_1+\ldots+\nu_d = 1 \quad\text{and}\quad \nu_1x_1+\ldots+\nu_dx_d= \frac1m\cdot v \vspace*{-1mm} \]
	and therefore $v = x_1^{\nu_1m}\cdot \ldots \cdot x_d^{\nu_dm}$.
	
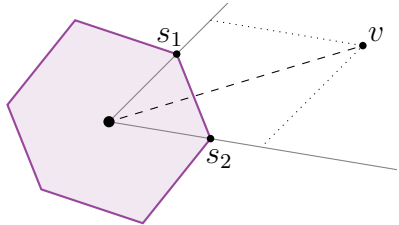
	\begin{center}
		\begin{tikzpicture}[scale=.45]
			\clip (-8.5,-3.1) rectangle (8.5,3.5);
			\draw[thick, Purple, fill=Purple!10] (3,-.5) -- (2,2) -- (-1,3) -- (-3,0.5) -- (-2,-2) -- (1,-3) -- cycle;
			\node (0) at (0,0) {};
			\node[circle, fill=black, inner sep=1pt] (v) at (7.5,2.25) {};
			\draw[dashed] (0) -- (v);
			\draw[gray] (12,-2) -- (0,0) -- (4,4);
			\draw[dotted] (3,3) -- (v) -- (4.5,-.75);
			\node[circle, fill=black, inner sep=1pt] at (2,2) {};
			\node[circle, fill=black, inner sep=1pt] at (3,-.5) {};
			\node[circle, fill=black, inner sep=1.5pt] at (0) {};
			\node at (7.9,2.6) {$v$};
			\node at (1.8,2.5) {$x_1$};
			\node at (3.25,-1.1) {$x_2$};
		\end{tikzpicture}
	\end{center}
	
	Next we show that any $\RR$-word $x_1^{\mu_1}\cdot \ldots \cdot x_k^{\mu_k}$ such that all letters $x_i$ lie on a common face $F$ are geodesics. Consider another $\RR$-word $y_1^{\lambda_1}\cdot\ldots\cdot y_\ell^{\lambda_\ell}\in X_\RR^*$ representing the same element as $x_1^{\mu_1}\cdot \ldots \cdot x_k^{\mu_k}$. As $F$ is a face, there exists a linear form $f\colon \RR^d \to \RR$ such that $f(p)\le 1$ for all $p\in P$, with equality if and only if $p\in F$. It follows that
	\[ \mu_1+\ldots+\mu_k = f(x_1^{\mu_1}\cdot \ldots \cdot x_k^{\mu_k}) = f(y_1^{\lambda_1}\cdot\ldots\cdot y_\ell^{\lambda_\ell}) = \lambda_1 f(y_1) + \ldots + \lambda_\ell f(y_\ell) \le \lambda_1+\ldots+\lambda_\ell \]
	which means $x_1^{\mu_1}\cdot \ldots \cdot x_k^{\mu_k}$ has indeed minimal length. In particular, this applies to the previously constructed $\RR$-word (of length $m$) for $v$: we have $\norm{v}_{\Stoll,X}=m$.
\end{proof}

\subsection{The $2$-step nilpotent case}

When $H$ is $2$-step nilpotent, the Malcev completion decomposes as $\bar H=V_1\oplus V_2$ and the group operation is given by
\[ (u,z_1)(v,z_2) = \Big(u+v,\; z_1+z_2 + \frac12[u,v]\Big), \]
where $[\cdot,\cdot]\colon V_1\times V_1\to V_2$ is a surjective, anti-symmetric bilinear map. Given an element $g\in \bar H$, we denote by $\pi(g)$ its first component and $z(g)$ the second component (its \say{areas}). Note that the areas of a product are given by
\begin{equation}\label{form:areas}
z(h_1h_2\ldots h_k) = \sum_{i=1}^k z(h_i) + \frac12 \sum_{i<j} \big[\pi(h_i),\pi(h_j)\big].
\end{equation}
In the $2$-step-nilpotent case, Stoll proves more than we could bargain for:
\begin{prop}[{Rough isometry, \cite[Proposition 4.3]{Stoll_metric}}] \label{prop:rough}
	Let $H$ be a torsion-free $2$-step nilpotent group, with a finite generating set $X$ and $\sigma\colon X\to \ZZ_{>0}$ a weight function. Then there exists a constant $C=C(X,\sigma)$ such that
	\[ \forall h\in H, \quad \norm{h}_{\Stoll,X,\sigma} \le \norm{h}_{X,\sigma} \le \norm{h}_{\Stoll,X,\sigma} + C. \]
\end{prop}
\medbreak
When considering words $w\in X^*$, there is usually a trade-off between having small coarse length $k(w)$, and having length $\ell_\sigma(w)$ close to $\norm{\overline w}_{X,\sigma}$. However, in the $2$-step nilpotent case, this trade-off does not happen. Analyzing the proofs of Lemma 4.2 and Proposition 4.3 of \cite{Stoll_metric}, we see that the word $w\in X^*$ evaluating to $h$, and witnessing $\norm{h}_{X,\sigma} \le \norm{h}_{\Stoll,X,\sigma} + C$, has bounded coarse length. More precisely,
\begin{lemma}[The best of both worlds] \label{sec1:lem_coarse}
	There exists a constant $K$ such that any element $h\in H$ can be written as $h=x_1^{m_1}x_2^{m_2}\ldots x_k^{m_k}$ with $x_i\in X$, $m_i\in \ZZ_{>0}$, $k\le K$ and 
	\[ \ell_\sigma(x_1^{m_1}x_2^{m_2}\ldots x_k^{m_k}) \le \norm{h}_{\Stoll,X,\sigma}+C \le \norm{h}_{X,\sigma} + C. \]
\end{lemma}
\textbf{Remark.} The Stoll paper only treats the case of symmetric generating set, with weight function $\sigma\equiv 1$. However the argument adapts to our more general setting.
\section{Generalities on virtually nilpotent groups} \label{sec:gen}

Let us start with some general observations on virtually nilpotent groups. In this section, $G$ is a group with a finite (non-weighted) generating set $S$. Let $H$ be a finite-index, torsion-free, $s$-step nilpotent subgroup. We do not require $H$ to be normal.

\subsection{A generating set for $H$} \label{sec:gen_set}

Consider the Schreier graph $\Sch(H\backslash G,S)$. We define a new generating (multi)set
\[ X(S) = \left\{ \,tat^{-1} \;\middle|\; \begin{array}{r} t\in S^* \text{ labels a simple path }H\;\to\, Ht \\ a\in S^* \text{ labels a simple cycle }Ht\to Ht \end{array} \right\}\subset H\]
with a weight function $\sigma\colon X\to\ZZ_{>0}$ defined by $\sigma(tat^{-1})=\ell(a)$.

\begin{lemma}[Decomposition map]
	There exists an injective map
	\[ \begin{pmatrix} \{w\in S^*\mid \bar w\in H\} & \longto & X^* \\ w & \longmapsto & \tilde w \end{pmatrix} \]
	such that $\overline w=\overline{\tilde w}$ and $\ell(w)=\ell_\sigma(\tilde w)$. In particular, $H=\la X\ra$.
\end{lemma}
\begin{proof} Consider $w\in S^*$ such that $\bar w\in H$. We describe its \emph{decomposition} $\tilde w \in X^*$ using a loop-erasure algorithm. We trace the path starting from $H$ and labeled by $w$ in $\Sch(H\backslash G,S)$, until we reach a vertex that was already visited, forming a loop $a_1$ that we erase, and then we keep going (forming other loops that we erase, ...).
\begin{center}
	\begin{tikzpicture}[scale=1.2]
	\draw[black, thick, -latex] (-1,-1)
		to [out=60, in=-180, looseness=.9] (0,0)
		to [out=0, in=20, looseness=2] (-.5,1)
		to [out=-160, in=-160, looseness=3] (-.8,1.5)
		to [out=20, in=100, looseness=2] (-.5,1)
		to [out=-80, in=120] (0,0)
		to [out=-60, in=170] (0.5,-.5)
		to [out=-10, in=-85, looseness=2] (1,-.2)
		to [out=95, in=90, looseness=2] (0.5,-.5)
		to [out=-90, in=0, looseness=1.5] (-1,-1);
	\node at (-.8,-.1) {$w$};
	\node at (1.4,0.3) {$=$};
	
	\draw[TealBlue, -latex, thick, shift={(3,0)}] (-1,-1)
		to [out=60, in=-180, looseness=.9] (0,0)
		to [out=0, in=20, looseness=2] (-.5,1);
	\draw[black, thick, -latex, shift={(3,0)}] (-.5,1)
		to [out=-160, in=-160, looseness=3] (-.8,1.5)
		to [out=20, in=100, looseness=2] (-.5,1);
	\node[TealBlue] at (2.2,-.1) {$t_1$};
	\node at (2.8,1.5) {$a_1$};
	\node at (3.7,0.3) {$\cdot$};
	
	\draw[TealBlue, -latex, thick, shift={(5,0)}] (-1,-1)
		to [out=60, in=-180, looseness=.9] (0,0);
	\draw[black, thick, -latex, shift={(5,0)}] (0,0)
		to [out=0, in=20, looseness=2] (-.5,1)
		to [out=-80, in=120] (0,0);
	\node[TealBlue] at (4.2,-.1) {$t_2$};
	\node at (5.1,1.2) {$a_2$};
	\node at (5.8,0.3) {$\cdot$};
	
	\draw[TealBlue, -latex, thick, shift={(7,0)}] (-1,-1)
		to [out=60, in=-180, looseness=.9] (0,0)
		to [out=-60, in=170] (0.5,-.5);
	\draw[black, thick, -latex, shift={(7,0)}] (0.5,-.5)
		to [out=-10, in=-85, looseness=2] (1,-.2)
		to [out=95, in=90, looseness=2] (0.5,-.5);
	\node[TealBlue] at (6.2,-.1) {$t_3$};
	\node at (7.5,0.2) {$a_3$};
	\node at (8.4,0.3) {$\cdot$};
	
	\draw[black, thick, -latex, shift={(9.7,0)}] (-1,-1)
		to [out=60, in=-180, looseness=.9] (0,0)
		to [out=-60, in=170] (0.5,-.5)
		to [out=-90, in=0, looseness=1.5] (-1,-1);
	\node at (9.7,.2) {$a_4$};
\end{tikzpicture} \vspace*{-3mm}
	\captionsetup{margin=10mm, font=footnotesize}
	
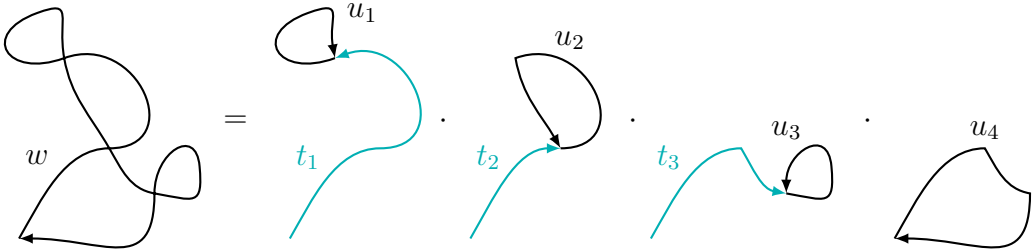
\captionof{figure}{A path decomposed as a product of \say{freeze frames} of the loop-easure algorithm.}
\end{center}
We call $t_i$ the word labelling the trace to $a_i$ in the process. We define
\[ \tilde w = (t_1a_1t_1^{-1}) \cdot (t_2a_2t_2^{-1}) \cdot \ldots \cdot (t_\ell a_\ell t_\ell^{-1}) \in X^* \]
We observe that $\overline w=\overline{\tilde w}$. (Note that $t_\ell$ is the empty word since $\bar w\in H$.) Moreover, since all the edges used by $w$ have been partitioned into simple cycles, $\ell(w)=\ell_\sigma(\tilde w)$.
	
Let us prove injectivity: we start with a decomposition $\tilde w\in X^*$ and reconstruct $w$. Observe that $t_i$ is a prefix of $t_{i+1}a_{i+1}$ for all $i<\ell$ and $t_\ell=\varepsilon$.\footnote{These two conditions characterize decompositions among words in $X^*$.} We define $r_i\in S^*$ such that $t_ir_i=t_{i+1}a_{i+1}$, then $w=t_1a_1r_1r_2\ldots r_{\ell-1}$. \vspace*{1mm}
\end{proof}

\subsection{Sub-linear control between word metrics} \label{sec:control}

Since $\ell_\sigma(\tilde w)=\ell(w)$ we have $\norm{h}_{X,\sigma}\le \norm{h}_S$ for all $h\in H$. The goal of this paragraph is to prove an inequality in the other direction. In order to state our result, we first define a sequence $(\alpha_s)_{s\ge 2}\subset [0,1)$. It starts with $\alpha_2=0$ and then
\[ \forall s\ge 3, \quad \alpha_s = \frac{1-\frac 1s\alpha_{s-1}}{2-\alpha_{s-1}-\frac 1s}. \]
An induction shows that $0\le\alpha_s<1$. We can now state our main inequality:
\begin{prop}\label{prop:sublin}
	Let $G$ be a virtually $s$-step nilpotent group, and $H$ a finite-index, torsion-free, $s$-step nilpotent subgroup of $G$. We consider $S$ a generating set for $G$ and $X$ the associated generating set for $H$, with the weight function $\sigma\colon X\to\ZZ_{>0}$. Then\footnote{When writing $O(g)$, we mean an expression $f$ satisfying $\abs{f}\le C g$ for a constant $C$ which only depends on $G,H,S,X,\sigma$. We often think of $O(g)$ as a positive error, hence write $-O(g)$ in some cases.}
	\[ \forall h\in H,\; \norm{h}_{X,\sigma} \le \norm{h}_S \le \norm{h}_{X,\sigma} + O\big(\!\norm{h}_{X,\sigma}^{\alpha_s}\big). \]
\end{prop}
We will need some preparatory results. First, recall the lower central series 
\[ \gamma_1(H)=H \quad\text{and, for each }i\ge 2,\quad  \gamma_i(H)=[\gamma_{i-1}(H),H]. \]
A classical lemma is the following, see for instance \cite[\S14.1.3]{GGT}
\begin{lemma}[Distortion]\label{lem:dist}
	Let $H$ be a finitely generated torsion-free $s$-step nilpotent group. Consider $\norm{\,\cdot\,}_E$ an Euclidean norm on $\gamma_s(H)\simeq\ZZ^c$. Then for $z\in\gamma_s(H)$
	\[ \norm z_{X,\sigma} = \Theta\big(\!\norm{z}_E^{1/s}\big) \quad\text{as }\norm z_E\to\infty. \]
\end{lemma}

Next we need the following generalization of Lemma \ref{sec1:lem_coarse}:
\begin{lemma}[$k$ versus $\ell_\sigma$]\label{lem:coarse_vs_length}
	Any element $h\in H$ with $\norm{h}_{X,\sigma}=n$ can be represented by a word $v\in X^*$ with coarse length $k(v) = O(n^{\alpha_s})$ and length $\ell_\sigma(v)= n + O(n^{\alpha_s})$.
\end{lemma}
\begin{proof}
	We argue by induction on $s$. The case $s=2$ (in which $\alpha_s=0$) is Lemma \ref{sec1:lem_coarse}.
	
	Suppose the induction hypothesis holds for $s-1\ge 2$. Consider $h\in H$ with $\norm h_{X,\sigma}=n$. Let $u\in X^*$ be a geodesic word representing $h$. We can decompose $u$ as a product $u = u_1u_2\ldots u_m$ with $m= n^\beta+O(1)$ pieces of length $n^{1-\beta}+O(1)$ for $\beta\in(0,1)$ which will be chosen later. By induction hypothesis, there exist words $v_i\in X^*$ such that
	\[ \overline v_i=\overline u_i \,\bmod{\gamma_s(H)}, \;\; k(v_i)= O\big(n^{(1-\beta)\alpha_{s-1}}\big)\;\;\text{and}\;\;  \ell_\sigma(v_i)=n^{1-\beta} + O\big(n^{(1-\beta)\alpha_{s-1}}\big).\]
	Observe that the error $z_i=\bar u_i\bar v_i^{-1}\in \gamma_s(H)$ has length 
	\[ \norm{z_i}_{X,\sigma}\le \ell_\sigma(u_i)+\ell_\sigma(v_i)=O(n^{1-\beta}), \vspace*{2mm} \]
	hence $\norm{z_i}_E=O(n^{(1-\beta)s})$ by Lemma \ref{lem:dist}. Therefore, the total error $z = z_1z_2\ldots z_m$ has size $\norm z_E=O(n^\beta\cdot n^{(1-\beta)s})$. The same lemma delivers $v_z\in X^*$ such that $\overline v_z=z$ and $\ell_\sigma(v_z)=O(n^{1-\frac{s-1}s\beta})$. Finally, let $v = v_1v_2\ldots v_m v_z \in X^*$. We have $\overline v=h$, and
	\begin{align*}
	k(v) & \le \sum_{i=1}^m k(v_i) + k(v_z) = n^\beta \cdot O(n^{(1-\beta)\alpha_{s-1}}) + O(n^{1-\frac{s-1}s\beta}), \\
	\ell_\sigma(v) & = \sum_{i=1}^m \ell_\sigma(v_i) +\ell_\sigma(v_z) = n + n^\beta \cdot O(n^{(1-\beta)\alpha_{s-1}}) + O(n^{1-\frac{s-1}s\beta}).	
	\end{align*}
	To conclude, we fine-tune $\beta=\frac{1-\alpha_{s-1}}{2-\alpha_{s-1}-\frac1s}$ so that $\beta+(1-\beta)\alpha_{s-1}=1-\frac{s-1}s\beta=\alpha_s$.
\end{proof}
\begin{rem}
	The first exponents $\alpha_2=0$ and $\alpha_3=\frac35$ are optimal (by Remark \ref{rem:BrLD}). Later $\alpha_s$ can probably be improved. For instance, Gianella has proved a result analogous to Lemma \ref{lem:coarse_vs_length} with $\alpha_s=\frac{s}{s+2}$ if we allow $w$ to be an $\RR$-word. \cite[Lemma 40]{Gianella}
\end{rem}
\begin{proof}[Proof of Proposition \ref{prop:sublin}]
	Let us prove the right inequality. Using the Lemma \ref*{lem:coarse_vs_length}, any element $h\in H$ of length $\norm h_{X,\sigma}=n$ can be represented by a word
	\[ v = (t_1a_1t_1^{-1})^{m_1} \dots (t_k a_k t_k^{-1})^{m_k} \in X^* \]
	with $k(v)=O(n^{\alpha_s})$ and $\ell_\sigma(v)=n+O(n^{\alpha_s})$. We convert this word into
	\[ w = t_1\,a_1^{m_1}\,u_1  \ldots t_k\,a_k^{m_k}\,u_k \in S^*, \]
	where $u_i\in S^*$ is a geodesic representative for $t_i^{-1}$. This word has length \vspace*{-3mm}
	\[
	\ell(w) = \ell_\sigma(v) + \sum_{i=1}^k \big(\ell(t_i)+\ell(u_i)\big) \le \ell_\sigma(v) + k(v) \big([G:H]+\max_{t}\|t^{-1}\|_S\big) = n+O(n^{\alpha_s}), \vspace*{-1mm}
	\]
	so that $\norm{h}_S\le \norm{h}_{X,\sigma}+O\big(\!\norm h_{X,\sigma}^{\alpha_s}\big)$ as announced.
\end{proof}

\subsection{A parte: Volume growth of virtually nilpotent groups}

Let us recall the state-of-the-art for growth of finitely generated nilpotent groups:
\begin{thm}[{\cite{Stoll_metric, Gianella}}] \label{thm:gianella}
	Let $H$ be a finitely generated $s$-step nilpotent group, and $X$ be a symmetric generating set. We have \vspace*{-2mm}
	\[ \beta_{(H,X)}(n) = c_{(H,X)} \cdot n^{d} + O(n^{d-\delta_s}), \vspace*{-3mm}\]
	where \vspace*{-1mm}
	\begin{itemize}[leftmargin=8mm]
		\item $d=d(H)=\sum_{i=1}^s i\cdot \mathrm{rank}_\QQ\big(\gamma_i(H)\big/\gamma_{i+1}(H)\big)\in\ZZ_{\ge0}$ is the Bass-Guivarc'h exponent,
		\item $c_{(H,X)}\in\RR_{>0}$ is the volume of the unit ball in the asymptotic cone of $(H,X)$ with its associated Carnot-Caratheodory metric.
		\item $\delta_s=1$ for $s=1,2$ and $\delta_s=\frac 1s$ for $s\ge 3$.
	\end{itemize}
\end{thm}
The statement extends if we add a weight function $\sigma\colon X\to\mathbb Z_{>0}$ into the picture. Our modest contribution is to extend this result to \emph{virtually} nilpotent groups
\begin{cor} \label{cor:vol}
	Let $G$ be a finitely generated, virtually $s$-step nilpotent group, and $S$ be a symmetric generating set. Let $H$ be a finite-index, torsion-free, $s$-step nilpotent group, with the associated generating set $X$ and weight function $\sigma\colon X\to\ZZ_{>0}$. We have
	\[ \beta_{(G,S)}(n) = [G:H]\cdot c_{(H,X,\sigma)} \cdot n^d + O(n^{d-\min\{1-\alpha_s;\;\delta_s\}}).\]
\end{cor}
\begin{proof}
	Let $j=[G:H]$, and decompose $G=\bigsqcup_{i=1}^j t_iH$. Picking $t_i$ as short as possible, we may assume that $\norm{t_i}_S<j$ for all $i$. We have
	\[ \bigsqcup_{i=1}^j \,t_i\cdot  \big(B_{(G,S)}(n-j) \cap H\big) \subseteq B_{(G,S)}(n) \subseteq \bigsqcup_{i=1}^j \,t_i\cdot  \big(B_{(G,S)}(n+j) \cap H\big) \]
	Combining this with Proposition \ref{prop:sublin}, we get the inclusions
	\[ \bigsqcup_{i=1}^j \;t_i\cdot B_{(H,X,\sigma)}\big(n-j-O(n^{\alpha_s})\big) \subseteq B_{(G,S)}(n) \subseteq \bigsqcup_{i=1}^j \;t_i\cdot B_{(H,X,\sigma)}(n+j) \]
	hence, by Theorem \ref{thm:gianella},
	\[ j \cdot c_{(H,X,\sigma)}\cdot \big(n-j-O(n^{\alpha_s})\big)^d + O(n^{d-\delta_s}) \le \beta_{(G,S)}(n) \le j \cdot c_{(H,X,\sigma)}\cdot (n+j)^d  + O(n^{d-\delta_s}), \]
	that is, $\beta_{(G,S)}(n)=[G:H]\cdot c_{(H,X,\sigma)}\cdot n^d + O(n^{d-\min\{1-\alpha_s;\;\delta_s\}})$.
\end{proof}

\subsection{Structure of almost-geodesics}

We give local conditions on almost geodesic $v\in X^*$. Essentially, most short subwords of $v$ should be geodesics in the abelianization $(\pi(\bar H),\norm{\,\cdot\,}_\Mink)$. More precisely
\begin{defi}
	Let $H$ be a torsion-free $s$-step nilpotent group, with a finite generating set $X$ and a weight function $\sigma\colon X\to \ZZ_{>0}$. Consider $\pi\colon H\onto H/[[H,H]]\simeq\ZZ^d$ and 
	\[ P=\Conv\left\{\frac{\pi(x)}{\sigma(x)}\;\middle|\; x\in X \right\}\subseteq\RR^d.\]
	A word $u\in X\cup X^2$ is called \emph{costly} if
	\begin{itemize}[leftmargin=8mm, label={\normalfont\textbullet}]
	\item $u=x$ with $\frac{\pi(x)}{\sigma(x)}$ not on the boundary of $P$, or
	\item $u=xy$ with $\frac{\pi(x)}{\sigma(x)},\frac{\pi(y)}{\sigma(y)}$ not on a common facet of $P$.
	\end{itemize}
	For $v\in X^*$, we define $N(v)$ as the number of occurrences of costly subwords in $v$.
\end{defi}
\begin{prop}\label{prop:blowup}
	There exists a constant $\delta=\delta(H,X,\sigma)>0$ such that, for all $v\in X^*$,
	\[ \ell_\sigma(v) \ge \norm{\overline v}_{X,\sigma} + \delta\cdot N(v) - O\big(\!\norm{\overline v}_{X,\sigma}^{\alpha_s}\big). \]
\end{prop}
In some sense, this result is a quantified, discrete analog to the \say{$(s-1)$-iterated blowup} result of Hakavuori and Le Donne. \cite[Corollary 1.4]{BlowUp}
\begin{proof}
	By Lemma \ref{sec1:abelian}, for each costly word $u$, there exists $u'\in X_\RR^*$ such that $\pi(u)=\pi(u')$ and that $\ell_\sigma(u)-\ell_\sigma(u')>0$. We fix a real number $\delta$ such that
	\[ 0<\delta < \frac12 \min\big\{ \ell_\sigma(u)-\ell_\sigma(u') \mid u\text{ is costly}\,\big\}. \]
	We now argue by induction on $s$. 
	
	$\blacktriangleright$ We first initialize for $s=2$: Consider a word $v\in X^*$ with $N=N(v)$ occurrences of costly subwords. Say $M\ge \frac12 N$ of these occurrences are disjoints, we denote them by $u_1,\ldots,u_M$. We replace $u_i$ in $v$ by $u_i'$, thus defining a new $\RR$-word $v'$.
	
	Observe that $v'$ has the same abelianization $\pi(\overline w')=\pi(\overline w)$. It only differs in areas:
	\[ z(\overline v')-z(\overline v) = \sum_{i=1}^M \big(z(\overline u_i')-z(\overline u_i)\big) \in [\bar H,\bar H]\simeq \RR^c \]
	(using Formula \ref{form:areas}). Recall that $[\bar H,\bar H]$ is quadratically distorted in $\bar H$ (Lemma \ref{lem:dist} but for Lie groups, see rather \cite[Lemme II.1]{Guivarc'h} or \cite[Theorem 2.7]{Breuillard}), hence 
	\[ \norm{z(\overline v')-z(\overline v)}_ {\Stoll,X,\sigma} = O\big(M^\frac12\big). \]
	Therefore, there exists an $\RR$-word $v_z$ with $\ell_\sigma(v_z)=O\big(M^\frac12)$ and $\overline v=\overline{v'v_z}$. It follows \vspace*{-2mm}
	\begin{IEEEeqnarray*}{rCLLCC}
	\ell_\sigma(v)
	& = &  \ell_\sigma(v'v_z) & + \sum_{i=1}^M \big(\ell_\sigma(u_i)-\ell_\sigma(u'_i)\big) & - & \ell_\sigma(v_z) \\
	& \ge & \norm{\overline v}_{\Stoll,X,\sigma} & + M \cdot \min\big\{ \ell_\sigma(u)-\ell_\sigma(u') \mid u\text{ is costly}\, \big\} & - & O\big(M^\frac12\big) \\[1.5mm]
	& \ge & \norm{\overline v}_{X,\sigma}-C & + \frac12N\cdot \min\big\{ \ell_\sigma(u)-\ell_\sigma(u') \mid u\text{ is costly}\, \big\} & - & O\big(N^\frac12\big) \\[1.5mm]
	& \ge & \norm{\overline v}_{X,\sigma} & + \delta\cdot N(v) \;-\; O(1),
	\end{IEEEeqnarray*}
	where we have used Proposition \ref{prop:rough} for the second inequality.
	
	$\blacktriangleright$ Suppose that the induction hypothesis holds for $s-1\ge 2$. Consider $v\in X^*$ of length $\ell_\sigma(v)=n$. We decompose $v$ as a product $v=v_1v_2\ldots v_m$ with $m= n^\beta+O(1)$ pieces of length $n^{1-\beta}+O(1)$. Using the hypothesis, there exists $u_i\in X^*$ such that
	\[ \overline u_i=\overline v_i \,\bmod{\gamma_s(H)} \quad\text{and}\quad  \ell_\sigma(u_i)\le \ell_\sigma(v_i) - \delta \cdot N(v_i) + O\big(n^{(1-\beta)\alpha_{s-1}}\big).\]
	We repeat part of the proof of Lemma \ref{lem:coarse_vs_length}: the error $z_i=\bar v_i\bar u_i^{-1}\in \gamma_s(H)$ has length 
	\[ \norm{z_i}_{X,\sigma}\le \ell_\sigma(u_i)+\ell_\sigma(v_i)=O(n^{1-\beta}) \]
	The same argument using distortion gives a word $u_z\in X^*$ such that $\overline u_z=z_1z_2\ldots z_m$ and $\ell_\sigma(u_z)=O(n^{1-\frac{s-1}s\beta})$. Finally, let $u = u_1u_2\ldots u_m u_z \in X^*$. We have $\overline u=\overline v$, and \vspace*{-1mm}
	\begin{align*}
	\norm{\overline v}_{X,\sigma}\le \ell_\sigma(u)
	& = \sum_{i=1}^m \ell_\sigma(u_i) + \ell_\sigma(u_z) \\
	& \le \sum_{i=1}^m\Big(\ell_\sigma(v_i)-\delta N(v_i)+O\big(n^{(1-\beta)\alpha_{s-1}}\big) \Big) + O(n^{1-\frac{s-1}s\beta}) \\
	& \le \ell_\sigma(v) - \delta\big(\vspace*{-.5mm}N(v)-m\big) + O\big(n^{\beta+(1-\beta)\alpha_{s-1}}\big) + O(n^{1-\frac{s-1}s\beta}).
	\end{align*}
	Fine-tuning $\beta=\frac{1-\alpha_{s-1}}{2-\alpha_{s-1}-\frac1s}$ gives us the desired result.
\end{proof}
\section{Proof of Theorem \ref{thm:crit}} \label{sec:main}
	From now on, we suppose $H$ is a torsion-free, $s$-step nilpotent, finite index, \emph{normal} subgroup of $G$.  As in Section \ref{sec:gen_set}, we consider the labeled graph $\Sch(H\backslash G,S)$ and
	\[ X(S) = \left\{ \,\bar t\bar a\bar t^{-1} \;\middle|\; \begin{array}{r} t\in S^* \text{ labels a simple path }H\;\to\, Ht \\ a\in S^* \text{ labels a simple cycle }Ht\to Ht \end{array} \right\}\subset H\]
	with a weight function $\sigma\colon X\to\ZZ_{>0}$ defined by $\sigma(\bar t\bar a\bar t^{-1})=\ell(a)$.
	\begin{rem*}
		Observe that $\Sch(H\backslash G,S)$ is transitive since $H$ is normal, hence a word $a\in S^*$ labels a simple cycle $Ht\to Ht$ if and only if it labels a simple cycle $H\to H$. We can therefore say a word $a\in S^*$ is a \emph{simple cycle} without ambiguity.
	\end{rem*}
	
	The proof naturally splits into two cases.
	
	$\blacktriangleright$ First, we suppose that two elements of $A(S)$ lie on the same facet of $P(S)$, and we conclude $(G,S)$ has exponential geodesic growth. We already provide exponentially many geodesics in the virtually abelian quotient $G/[[H,H]]$. Observe that
	\[ P(S) \overset{\mathrm{def}}= \Conv\left(\left\{ \frac{\pi(\bar a)}{\ell(a)} \;\middle|\; a\in S^* \text{ is a simple cycle} \right\}^{G/H}\right) = \Conv\left\{ \frac{\pi(x)}{\sigma(x)} \;\middle|\; x\in X\right\}.\]
	is the polytope that governs the Minkowski norm on $\hat H=H/[[H,H]]\simeq \ZZ^d$ with respect to the natural generating multiset $\hat X=\pi(X)$, with the weight function $\hat\sigma(\hat x)=\sigma(x)$. 
	\medbreak
	Consider distinct simple cycles $a,b\in S^*$ (in particular $\bar a,\bar b\in X$) such that  \vspace*{-1mm}
	\[ \frac{\pi(\bar a)}{\ell(a)}= \frac{\pi(\bar a)}{\sigma(a)} \quad\text{and}\quad \frac{\pi(\bar b)}{\ell(b)} = \frac{\pi(\bar b)}{\sigma(b)} \]
	lie on a common facet of $P(S)$. Then, for all $w\in\{a,b\}^*$, we have
	\[ \ell(w) = \ell_{\hat \sigma}(\hat w) = \norm{\pi(\bar w)}_{\Mink,P} \le \norm{\pi(\bar w)}_{\hat X,\hat\sigma}\le \norm{\bar w}_{X,\sigma} \le \norm{\bar w}_S.\]
	(The second equality is justified by Lemma \ref{sec1:abelian}, and the last inequality by the easy part of Proposition \ref{prop:sublin}.) It follows that all those words are geodesics, which concludes. \hfill$\blacktriangleleft$
	
	\bigbreak
	
	From now on, we work towards an upper bound. First, we work specifically on geodesics $w\in S^*$ evaluating in the subgroup $H$.
	
	$\blacktriangleright$ Under the hypothesis that no two elements of $A(S)$ lie on the same facet of $P(S)$, we prove that the coarse length of the decomposition $\tilde w\in X^*$ (defined in \S\ref{sec:gen_set}) satisfies
	\begin{equation} \label{eq:changes_are_costly}
		k(\tilde w) \le N(\tilde w) \cdot [G:H]+1
	\end{equation}
	for any word $w\in S^*$ evaluating in $H$. Let us write
	\[ \tilde w = x_1^{m_1} \, x_2^{m_2}\cdot \ldots \cdot x_k^{m_k} \quad\text{with }x_i\ne x_{i+1}\text{ for all }i. \]
	The decomposition process not only gives this sequence of generators $(x_i)\in X$, but also two sequences of simple paths $(t_i)\in S^*$ and simple cycles $(a_i)\in S^*$ such that $x_i=\bar t_i\bar a_i\bar t_i^{-1}$. Observe that, for each $i$, one of $t_i$ and $t_{i+1}$ is a prefix of the other (depending on where the walk re-intersects itself). Consider a time $i$ such that $t_{i+1}$ is a prefix of $t_i$ (including the case $t_{i+1}=t_i$). In particular, we can rewrite $x_{i+1} = \bar t_i\bar b_{i+1}\bar t_i^{-1}$ for some simple cycle $b_{i+1}\in S^*$. ($b_{i+1}$ is a cyclic permutation of $a_{i+1}$.)
	\begin{center}
		\begin{tikzpicture}[scale=1.4]
\begin{scope}[shift={(0,0)}]
	\draw[Gray, very thick] (-.2,-1.4)
		to [out=45, in=-110, looseness=1.5] (0,0)
		to [out=160, in=-20, looseness=1.5] (-1,.1)
		to [out=-100, in=90, looseness=1.5] (-.67,-.7);
	\draw[TealBlue, very thick, -latex] (-.67,-.7)
		to [out=-90, in=-20, looseness=2] (-.95,-.95)
		to [out=160, in=-165, looseness=2] (-.67,-.7);
	\draw[Gray, dotted] (0,0)
		to [out=70, in=150, looseness=1.2] (.5,.5)
		to [out=-30, in=-20, looseness=1.2] (0,0);
	\draw[Gray, dotted] (-1,.1)
		to [out=160, in=-135, looseness=1.5] (-1.1,.5)
		to [out=45, in=80, looseness=1.5] (-1,.1);
	\node at (.25,-.95) {$t_i$};
	\node at (-1.25,-1) {$a_i$};
\end{scope}	

\begin{scope}[shift={(2.3,0)}]
	\draw[Gray, very thick] (-.2,-1.4)
		to [out=45, in=-110, looseness=1.5] (0,0);
	\draw[Purple!70, very thick, -latex] (0,0)
		to [out=160, in=-20, looseness=1.5] (-1,.1)
		to [out=-100, in=90, looseness=1.5] (-.67,-.7)
		to [out=25, in=-145, looseness=2] (-.03,-.012);
	\draw[Gray, dotted] (-.67,-.7)
		to [out=-90, in=-20, looseness=2] (-.95,-.95)
		to [out=160, in=-165, looseness=2] (-.67,-.7);
	\draw[Gray, dotted] (0,0)
		to [out=70, in=150, looseness=1.2] (.5,.5)
		to [out=-30, in=-20, looseness=1.2] (0,0);
	\draw[Gray, dotted] (-1,.1)
		to [out=160, in=-135, looseness=1.5] (-1.1,.5)
		to [out=45, in=80, looseness=1.5] (-1,.1);
	\node at (.4,-.95) {$t_{i+1}$};
	\node at (-1.1,-.5) {$a_{i+1}$};
\end{scope}

\begin{scope}[shift={(4.6,0)}]
\draw[Gray, very thick] (-.2,-1.4)
	to [out=45, in=-110, looseness=1.5] (0,0)
	to [out=160, in=-20, looseness=1.5] (-1,.1)
	to [out=-100, in=90, looseness=1.5] (-.67,-.7);
\draw[Orange, very thick, -{Latex[length=2.2mm, width=2mm, left]}] (-.67,-.7)
	to [out=25, in=-145, looseness=1.9] (-.05,-.02)
	to [out=161, in=-20, looseness=1.4] (-0.97,.052)
	to [out=-97, in=100, looseness=1.15] (-.63,-.68);	
\draw[Gray, dotted] (-.67,-.7)
	to [out=-90, in=-20, looseness=2] (-.95,-.95)
	to [out=160, in=-165, looseness=2] (-.67,-.7);
\draw[Gray, dotted] (0,0)
	to [out=70, in=150, looseness=1.2] (.5,.5)
	to [out=-30, in=-20, looseness=1.2] (0,0);
\draw[Gray, dotted] (-1,.1)
	to [out=160, in=-135, looseness=1.5] (-1.1,.5)
	to [out=45, in=80, looseness=1.5] (-1,.1);
\node at (.25,-.95) {$t_i$};
\node at (-1.1,-.5) {$b_{i+1}$};
\end{scope}		
\end{tikzpicture}
		\hspace*{5mm}
		\begin{tikzpicture}[scale=1.4]
	\draw[Gray, very thick] (-.2,-1.4)
		to [out=45, in=-110, looseness=1.5] (-.02,-.04);
	\draw[TealBlue, very thick] (.02,.04)
		to [out=70, in=150, looseness=1.2] (.5,.5)
		to [out=-30, in=-20, looseness=1.2] (0,0);
	\draw[Purple!70, very thick, -latex] (0,0)
		to [out=160, in=80, looseness=1.5] (-.6,-.2)
		to [out=-100, in=-145, looseness=1.5] (-.04,-.02);
	
	\node at (0.6,-.9) {\footnotesize$t_i=t_{i+1}$};
	\node at (.5,-.15) {$a_i$};
	\node at (-.85,-.55) {\footnotesize$a_{i+1}=b_{i+1}$};
\end{tikzpicture}
		\captionsetup{margin=4mm, font=footnotesize}
		
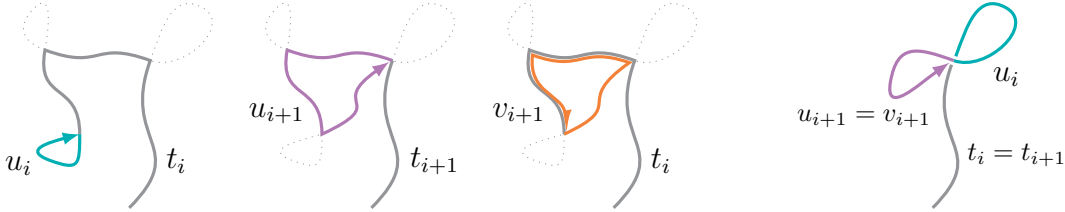
\captionof{figure}{$x_i=t_i\,a_i\,t_i^{-1}$ and $x_{i+1}=t_{i+1}\,a_{i+1}\,t_{i+1}^{-1}=t_i\,b_{i+1}\,t_i^{-1}$. \hspace*{7mm} The limit case $t_i=t_{i+1}$.}
	\end{center}
	As $x_i\ne x_{i+1}$, we have $a_i\ne b_{i+1}$. Now our hypothesis on $A(S)$ kicks in: both points
	\[ \frac{\pi(\bar a_i)}{\ell(a_i)} \quad\text{and}\quad \frac{\pi(\bar b_{i+1})}{\ell(b_{i+1})} \]
	cannot lie on a common facet of $P$. We get the same conclusion for 
	\[ \frac{\pi(x_i)}{\sigma(x_i)} \quad\text{and}\quad \frac{\pi(x_{i+1})}{\sigma(x_{i+1})} \]
	after conjugation by $t_i$. We conclude that the subword $x_ix_{i+1}$ is \emph{costly} in the sense of Proposition \ref{prop:blowup} (hence counted in $N(\tilde w)$), as soon as $t_{i+1}$ is a prefix of $t_i$.
	
	In order to avoid this situation, the path $t_i$ should be a proper prefix of $t_{i+1}$. This means the $t_i$'s usually get longer and longer, but they have lengths bounded between $0$ and $[G:H]-1$. Combining these two observations, we get the desired bound.
	
	\medbreak
	
	$\blacktriangleright$ We are now able to combine everything. Consider a word $w\in S^*$ representing an element $h\in H$. Propositions \ref{prop:sublin} and \ref{prop:blowup} give us
	\begin{align*}
	\norm{h}_S & \le \norm{h}_{X,\sigma} + O\big(\!\norm{h}_{X,\sigma}^{\alpha_s}\big) \\
	\ell(w) & = \ell_\sigma(\tilde w) \ge \norm{h}_{X,\sigma} + \delta\cdot N(\tilde w) - O\big(\! \norm{h}_{X,\sigma}^{\alpha_s}\big)
	\end{align*}
	for some constant $\delta>0$. If $w$ is a geodesic, then $\ell(w)=\norm{h}_S$ and therefore
	\[ \delta\cdot N(\tilde w) = O\big(\!\norm{h}_{X,\sigma}^{\alpha_s}\big) = O\big(\ell(\tilde w)^{\alpha_s}\big). \]
	Combined with the inequality (\ref{eq:changes_are_costly}), we get a constant $C=C(G,S)>0$ such that
	\begin{equation} \label{eq:geodesic_have_few_changes}
		k(\tilde w)\le C\cdot \ell(\tilde w)^{\alpha_s}.
	\end{equation}
	
	\medbreak
$\blacktriangleright$ We obtain the desired upper bound. The inequality (\ref*{eq:geodesic_have_few_changes}) gives an injection
\[ \left\{w\in S^* \;\big|\; w \text{ is geodesic with }\overline w\in H\right\} \into \left\{v\in X^* \;\big|\; k(v)\le C\cdot \ell(v)^{\alpha_s} \right\} \hspace*{3.5mm} \]
sending $w\mapsto\tilde w$. Observe that $\ell(\tilde w)\le \ell_\sigma(\tilde w)=\ell(w)$, hence
\[ \#\!\left\{w\in S^* \,\big|\, w \text{ is geodesic, }\ell(w)\le n,\; \overline w\in H\right\} 
\le  \#\!\left\{v\in X^* \,\big|\, k(v)\le C\cdot n^{\alpha_s},\;\ell(v)\le n\right\}. \]
Note that, in to order to construct a word $v$ satisfying these conditions, it suffices to pick $Cn^{\alpha_s}$ cut points along an interval of length $n$, pick a letter to fill each of the first $Cn^{\alpha_s}$ \say{blocks}, and leave the last block empty. \vspace*{-2mm}
\begin{center}
	\begin{tikzpicture}[scale=.8]
		\foreach \x in {1,...,12}{
		\draw[thick, rounded corners=.1] ({\x-.35},0) -- ({\x-.35},-.2) -- ({\x+.35},-.2) -- ({\x+.35},0);
		}
		\node at (.5,0) {$|$};
		\node at (3.45,0) {$|$};
		\node at (3.55,0) {$|$};
		\node at (6.5,0) {$|$};
		\node at (7.5,0) {$|$};
		\node at (10.5,0) {$|$};
		\node at (1,0) {$x_3$};
		\node at (2,0) {$x_3$};
		\node at (3,0) {$x_3$};
		\node at (4,0) {$x_1$};
		\node at (5,0) {$x_1$};
		\node at (6,0) {$x_1$};
		\node at (7,0) {$x_1$};
		\node at (8,0) {$x_4$};
		\node at (9,0) {$x_4$};
		\node at (10,0) {$x_4$};
	\end{tikzpicture}\vspace*{-2mm}
\end{center}
This translates into a crude upper bound of
\[
\#\!\left\{v\in X^* \,\big|\, k(v)\le C\cdot n^{\alpha_s},\;\ell(v)\le n\right\}
\le \big((n+1)\abs X\big)^{Cn^{\alpha_s}}
\asymp \exp\!\big(Cn^{\alpha_s}\log(n)\big).
\]
To conclude, each geodesic $w$ in $(G,S)$ can be decomposed as a product
\[ w = w_1s_1w_2s_2\ldots s_{k-1}w_k, \]
where each $w_i\in S^*$ is a geodesic ending up in $H$, $s_i\in S$, and $k\le [G:H]$, which gives an upper bound on $\gamma_\geod(n)$ of the same type.
\hfill$\qedsymbol$
\section{Geometry of the Engel group} \label{sec:Eng}

In this section, we are looking at a $3$-step nilpotent group, the Engel group $\Ec$. We first re-introduce a nilpotent Lie group $\bar\Ec$ using a geometrical model. $\Ec$ is defined as particular a cocompact lattice in $\bar\Ec$.

We consider $\Gamma$ the monoid of absolutely continuous paths in $\RR^2$ starting from $(0,0)$, with concatenation as operation. Then $\bar\Ec$ is defined as the quotient $\Gamma/\sim$ for an equivalence relation $\sim$. For any path $g\in\Gamma$, we define three parameters:
\begin{itemize}[leftmargin=8mm, label=\textbullet]
	\item[(1)] its second endpoint $\hat g=(x_g,y_g)\in\RR^2$.
	\item[(\textbullet)] a distribution of winding numbers. First, we get a closed path $g_c$ by concatenating $g$ with the segment back from $\hat g$ to $(0,0)$. Then the function $W_g\colon\RR^2\setminus \mathrm{Im}(g_c)\to \ZZ$ is defined as $W_g(x,y)=$ the winding number of $g_c$ around $(x,y)$. (See Figure \ref{fig:Engel_def}.)
	\item[(2)] its total algebraic (or signed) area \vspace*{1mm}
	\[ A(g) = \displaystyle \iint_{\RR^2} W_g(x,y)\;\mathrm dx\,\mathrm dy. \vspace*{-1mm} \]
	\item[(3)] the $y$-coordinate of its \say{barycenter} (or center of gravity)
	\[ B_y(g) = \displaystyle\iint_{\RR^2} y\cdot W_g(x,y) \;\mathrm dx\,\mathrm dy. \]
\end{itemize}
Two paths $g,h$ are equivalent if they share same endpoint $\hat g=\hat h$, same algebraic area $A(g)=A(h)$ and same \say{$y$-coordinate of barycenter} $B_y(g)=B_y(h)$.

\begin{prop}Given two paths $g,h$, their concatenation $gh$ has parameters
	\begin{align*}
	\widehat{gh} \hspace*{1.5mm} & = \widehat g + \widehat h, \\
	A(gh) & = A(g) + A(h) + \frac12\det\!\big(\hat g,\hat h\big), \\
	B_y(gh) & = B_y(g)+ B_y(h) + y_g\cdot A(h) + \frac13(2y_g+ y_h)\cdot\frac12\det\!\big(\hat g,\hat h\big).
	\end{align*}
	As a corollary, the operation \say{concatenation} passes to the quotient. With the empty path as neutral element and reverse paths as inverses, this defines a group.
\end{prop}
\begin{proof} The relation $\widehat{gh}=\hat g+\hat h$ is obvious. The key observation is the decomposition
	\[ W_{gh} = W_g + W_h\circ \tau_{-\hat g} \pm \chi_{\triangle((0,0),\hat g,\hat g+\hat h)}, \]
	where
	\begin{itemize}
		\item $\tau_v\colon\RR^2\to\RR^2$ is the translation by $v$,
		\item $\triangle\big((0,0),\hat g,\hat g+\hat h\big)$ is the convex hull of those three points, $\chi_{\triangle((0,0),\hat g,\hat g+\hat h)}$ denotes its characteristic function, and the sign $\pm$ depends on the order of $(0,0)$, $\hat g$ and $\hat g+\hat h$ on the boundary of the triangle ($-1$ if clockwise and $+1$ if anti-clockwise).
	\end{itemize}
	Here is a pictorial explanation: \vspace*{-4mm}
	\begin{center}
		\begin{tikzpicture}[scale=1.32]
			\clip (-.8,-.5) rectangle (2.55,3.4);
			\fill[Purple!40, opacity=.5] (0,0)
				to[out=-10, in=-110, looseness=1.5] (1.2,.6)
				to[out=70, in=160, looseness=2] (2.2,.8);
			\draw[Purple, thick, -latex] (0,0)
				to[out=-10, in=-110, looseness=1.5] (1.2,.6)
				to[out=70, in=160, looseness=2] (2.2,.8);
			\draw[dashed, Purple, thick, latex-] (0,0) -- (2.2,.8);
			
			\node[Purple] at (.7,.11) {\scriptsize$+1$};
			\node[Purple] at (1.55,.8) {\scriptsize$-1$};
		
			\fill[orange!40, opacity=.5] (0,0)
				to [out=70, in=-60, looseness=1.5] (0,1.3)
				to [out=120, in=65, looseness=1.5] (-.4,1.2)
				to [out=-115, in=-90, looseness=1.2] (0,1.3)
				to [out=90, in=65, looseness=2] (-.5,1.3)
				to [out=-115, in=-150, looseness=2] (.2,.9)
				to [out=30, in=-90, looseness=1.5] (1,2);
			\fill[orange!40]	(0,1.3)
				to [out=120, in=65, looseness=1.5] (-.4,1.2)
				to [out=-115, in=-90, looseness=1.2] (0,1.3);
			\draw[orange, thick, -latex] (0,0)
				to [out=70, in=-60, looseness=1.5] (0,1.3)
				to [out=120, in=65, looseness=1.5] (-.4,1.2)
				to [out=-115, in=-90, looseness=1.2] (0,1.3)
				to [out=90, in=65, looseness=2] (-.5,1.3)
				to [out=-115, in=-150, looseness=2] (.2,.9)
				to [out=30, in=-90, looseness=1.5] (1,2);
			\draw[dashed, orange, thick, latex-] (0,0) -- (1,2);
			
			\node[Orange] at (-.2,.96) {\scriptsize$+1$};
			\node[Orange] at (-.2,1.25) {\scriptsize$+2$};
			
			\node[circle, fill=black, inner sep=1.5pt, label=below:{\footnotesize$(0,0)$}] at (0,0) {};
			\node[circle, fill=Purple, inner sep=1pt, label=below:{\color{Purple}$\hat g$}] at (2.2,.8) {};
			\node[circle, fill=orange, inner sep=1pt, label=above:{\color{orange}$\hat h$}] at (1,2) {};
		\end{tikzpicture} \hspace*{2mm}
		\begin{tikzpicture}[scale=1.32]
			\clip (-.35,-.5) rectangle (3.55,3.4);
			
			\fill[TealBlue!40, opacity=.5] (0,0) -- (2.2,.8) -- (3.2,2.8) -- cycle;
			\draw[dashed, TealBlue, thick, latex-] (0,0) -- (3.2,2.8);
			
			\fill[Purple!40, opacity=.5] (0,0)
				to[out=-10, in=-110, looseness=1.5] (1.2,.6)
				to[out=70, in=160, looseness=2] (2.2,.8);
			\draw[Purple, thick, -latex] (0,0)
				to[out=-10, in=-110, looseness=1.5] (1.2,.6)
				to[out=70, in=160, looseness=2] (2.2,.8);
			\draw[dashed, Purple, thick] (0,0) -- (2.2,.8);
			\node[Purple] at (.7,.11) {\scriptsize$+1$};
			\node[Purple] at (1.55,.8) {\scriptsize$-1$};
			
			\begin{scope}[shift={(2.2,.8)}]
				\fill[orange!40, opacity=.5] (0,0)
					to [out=70, in=-60, looseness=1.5] (0,1.3)
					to [out=120, in=65, looseness=1.5] (-.4,1.2)
					to [out=-115, in=-90, looseness=1.2] (0,1.3)
					to [out=90, in=65, looseness=2] (-.5,1.3)
					to [out=-115, in=-150, looseness=2] (.2,.9)
					to [out=30, in=-90, looseness=1.5] (1,2);
				\fill[orange!40]	(0,1.3)
					to [out=120, in=65, looseness=1.5] (-.4,1.2)
					to [out=-115, in=-90, looseness=1.2] (0,1.3);
				\draw[orange, thick, -latex] (0,0)
					to [out=70, in=-60, looseness=1.5] (0,1.3)
					to [out=120, in=65, looseness=1.5] (-.4,1.2)
					to [out=-115, in=-90, looseness=1.2] (0,1.3)
					to [out=90, in=65, looseness=2] (-.5,1.3)
					to [out=-115, in=-150, looseness=2] (.2,.9)
					to [out=30, in=-90, looseness=1.5] (1,2);
				\draw[dashed, orange, thick] (0,0) -- (1,2);
				\node[Orange] at (-.2,.96) {\scriptsize$+1$};
				\node[Orange] at (-.2,1.25) {\scriptsize$+2$};
			\end{scope}
			
			\node[circle, fill=black, inner sep=1.5pt, label=below:{\footnotesize$(0,0)$}] at (0,0) {};
			\node[circle, fill=Purple, inner sep=1pt, label=below:{\color{Purple}$\hat g$}] at (2.2,.8) {};
			\node[circle, fill=orange, inner sep=1pt, label=above:{\footnotesize$\hat g+\hat h$}] at (3.2,2.8) {};
			
			\node[TealBlue] at (2,1.25) {\footnotesize$+1$};
		\end{tikzpicture}
		\begin{tikzpicture}[scale=1.32]
		\clip (-.35,-.5) rectangle (3.55,3.4);
		
		\begin{scope}[shift={(2.2,.8)}]
		\fill[TealBlue!40, opacity=.5] (-2.2,-.8)
			to[out=-10, in=-110, looseness=1.5] (-1,-.2)
			to[out=70, in=160, looseness=2] (0,0)
			to [out=70, in=-60, looseness=1.5] (0,1.3)
			to [out=120, in=65, looseness=1.5] (-.4,1.2)
			to [out=-115, in=-90, looseness=1.2] (0,1.3)
			to [out=90, in=65, looseness=2] (-.5,1.3)
			to [out=-115, in=-150, looseness=2] (.2,.9)
			to [out=30, in=-90, looseness=1.5] (1,2);
		\fill[TealBlue!40]	(0,1.3)
			to [out=120, in=65, looseness=1.5] (-.4,1.2)
			to [out=-115, in=-90, looseness=1.2] (0,1.3);
		\begin{scope}
			\clip (-2.2,-.8) -- (1,2) -- (.2,.9) -- cycle;
			\fill[TealBlue!40] (0,0)
				to [out=70, in=-60, looseness=1.5] (0,1.3)
				to [out=120, in=65, looseness=1.5] (-.4,1.2)
				to [out=-115, in=-90, looseness=1.2] (0,1.3)
				to [out=90, in=65, looseness=2] (-.5,1.3)
				to [out=-115, in=-150, looseness=2] (.2,.9);
		\end{scope}
		\draw[TealBlue, thick, -latex] (-2.2,-.8)
			to[out=-10, in=-110, looseness=1.5] (-1,-.2)
			to[out=70, in=160, looseness=2] (0,0)
			to [out=70, in=-60, looseness=1.5] (0,1.3)
			to [out=120, in=65, looseness=1.5] (-.4,1.2)
			to [out=-115, in=-90, looseness=1.2] (0,1.3)
			to [out=90, in=65, looseness=2] (-.5,1.3)
			to [out=-115, in=-150, looseness=2] (.2,.9)
			to [out=30, in=-90, looseness=1.5] (1,2);
			\node[TealBlue] at (-.2,1.25) {\scriptsize$+2$};
			\node[TealBlue] at (-.02,.94) {\scriptsize$+2$};
			\node[TealBlue] at (-.41,.98) {\scriptsize$+1$};
		\end{scope}
		\node[TealBlue] at (.7,.11) {\scriptsize$+1$};
		\node[TealBlue] at (1.55,.8) {\scriptsize$0$};
		
		\draw[dashed, TealBlue, thick, latex-] (0,0) -- (3.2,2.8);
		\node[circle, fill=black, inner sep=1.5pt, label=below:{\footnotesize$(0,0)$}] at (0,0) {};
		\node[circle, fill=TealBlue, inner sep=1pt, label=above:{\footnotesize$\hat g+\hat h$}] at (3.2,2.8) {};
		
		\node[TealBlue] at (2,1.25) {\footnotesize$+1$};
		\end{tikzpicture}
		
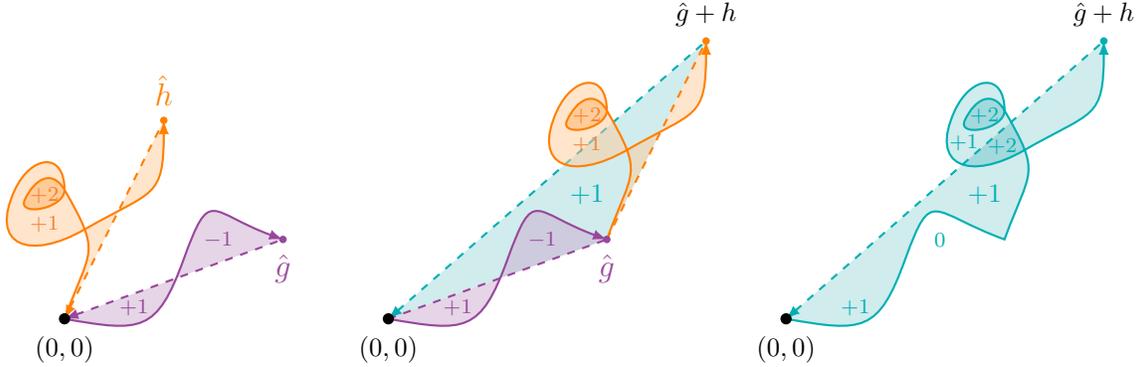
\captionof{figure}{Pictures of $g$, $h$ and $gh$ and some winding numbers.} \label{fig:Engel_def}
	\end{center}
It only remains to compute $A(gh)$ and $B_y(gh)$ using the decomposition. Observe that $\frac12\det(\hat g,\hat h)$ is the signed area of the triangle $\triangle\big((0,0),\hat g,\hat g+\hat h\big)$, and the $y$-coordinate of its center of gravity is given by $\frac13\big(0+y_g+(y_g+y_h)\big)$.
\end{proof}
\bigbreak
As our main example, we will consider the lattice generated by the straight segments $a$ and $b$ from $(0,0)$ to $(1,1)$ and $(1,-1)$ respectively. This group is given by \vspace*{-1mm}
\[ \Ec = \la a,b\;\big|\; [a,[a,b]]=[b^{-1},[a,b]]\text{ is central}\hspace*{2pt} \ra. \]
Throughout, we fix $X=\{a^\pm,b^\pm\}$ as a generating set. \vspace*{-2mm}
\begin{center}
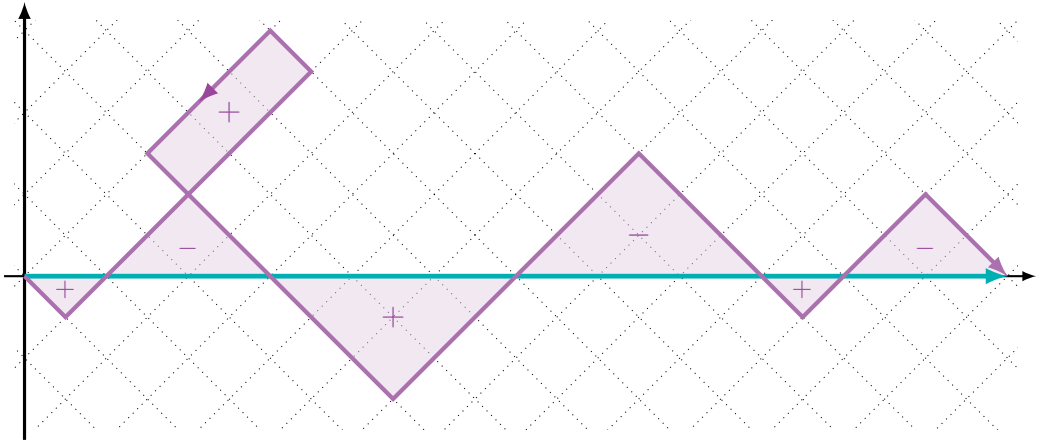

	\begin{tikzpicture}[scale=.72, rotate=-45]
		\draw[very thick, -latex] (2,-2) -- (-3.35,3.35);
		\draw[thick, -latex] (-.25,-.25) -- (12.35,12.35);
		
		\clip (1.75,-2) -- (14,10.25) -- (9,15.25) -- (-3.25,3) -- cycle;
		\draw[thin, dotted] (-4,-2) grid (15,16);
		
		\draw[Purple, fill=Purple!20, fill opacity=.5] (0,0) -- (1,0) -- (1,6) -- (0,6) -- (0,3) -- (6,3) -- (6,9) -- (10,9) -- (10,12) -- (12,12);
		\draw[line width=1.8pt, TealBlue, -latex] (0,0) -- (12,12);
		\draw[ultra thick, Purple!75, -latex] (0,0) -- (1,0) -- (1,6) -- (0,6) -- (0,3) -- (6,3) -- (6,9) -- (10,9) -- (10,12) -- (12,12);
		\draw[ultra thick, Purple, -latex] (0,4.5) -- (0,4.25);
	
		\node[Purple] at (.66,.33) {\footnotesize$\bf +$};
		\node[Purple] at (1.66,2.33) {\footnotesize$\bf -$};
		\node[Purple] at (.5,4.5) {$\bf +$};
		\node[Purple] at (5,4) {$\bf +$};
		\node[Purple] at (7,8) {$\bf -$};
		\node[Purple] at (9.66,9.33) {\footnotesize$\bf +$};
		\node[Purple] at (10.66,11.33) {\footnotesize$\bf -$};
	\end{tikzpicture}\vspace*{-2mm}
	\captionsetup{margin=5mm, font=footnotesize}
	\captionof{figure}{A word over $X$ and an equivalent path. Both satisfy $\hat g=(24,0)$ and $A(g)=B_y(g)=0$.} \vspace*{1mm}
\end{center}

\subsection{An observation of Stoll and lower bound on word length}

In unpublished notes \cite{Stoll_3step}, Stoll shows that a key result (\cite [Lemma 3.3]{Stoll_metric}: \say{Every element of $\bar H$ is represented by an $\RR$-word of minimal length}, whenever $\bar H$ is a simply connected $2$-step nilpotent group) fails in the Engel group\footnote{Stoll works instead with the Cartan group, i.e., the free $3$-step nilpotent group of rank $2$, but his argument works just as well for the Engel group, and I would rather not introduce new notations.} $\bar\Ec$. Specifically

\begin{prop*}[\cite{Stoll_3step}]
The element $g_1$ represented by a segment from $(0,0)$ to $(1,0)$ has Stoll length $\norm{g_1}_{\Stoll,X}=1$, but is not represented by any $\RR$-word of length $1$.
\end{prop*}
By a compactness argument, this proves that Lemma \ref{sec1:lem_coarse} does not extend either:

\begin{cor*}
Consider $g_n\in\overline\Ec$ the element represented by a segment to $(n,0)$. There does not exist a sequence of $\RR$-words $w_n$ representing $g_n$ such that the coarse lengths $k(w_n)$ are uniformly bounded and $\ell(w_n)\le\norm{g_n}_{\Stoll,X}+o(n)=n+o(n)$.	
\end{cor*}


\begin{proof}[Proof of the Corollary]
	We argue by contradiction, and suppose such a sequence of words does exist. After rescaling by a factor of $\frac1n$, we get a sequence of $\RR$-words of length $1+o(1)$ and coarse length $\le K$, all representing $g_1$. However, $\RR$-words with length $\le 2$ and coarse length $\le K$ form a compact set, so some subsequence converges to an $\RR$-word representing $g_1$ of length $1$ and coarse length $\le K$, a contradiction.
\end{proof}

In what follows, we quantify the dependence between $k(w_n)$ and $\ell(w_n)-\norm{g_n}_\Stoll$. We use that the horizontal path is a (particularly bad) abnormal curve in $\overline\Ec$.\vspace*{1mm}


\begin{prop}[$B_y$ back from the depths]\label{prop:back_from_depths} There exists $C>0$ such that, for any $w\in X_\RR^*$ with endpoint $\hat w=(n,0)$ and length $\ell(w)=n+\Delta$, we have \vspace*{-.5mm}
\[ -B_y(\bar w) \ge \frac1{24} \frac{n^3}{\big(k(w)-1\big)^2} - C\cdot \Delta^2 \cdot \max\left\{\Delta;\, \frac n{k(w)-1}\right\}.\]
\end{prop}
\begin{proof} We decompose the path $w$ into a main path (in purple), and some loops and boundary mess (in orange), and estimate the contribution of each part to $-B_y(\bar w)$. \vspace*{-7mm}
\begin{center}
\begin{tikzpicture}[scale=1.2]
	\clip (-1.5,-5.2) rectangle (11.8,2.2);
	
	\draw[very thick, -latex] (-1.1,0) -- (11,0);
	
	\draw[magenta, very thick] (.4,0) -- (.8,.4) -- (1.2,0) -- (.3,-.9) -- (-.6,0) -- (.4,1) -- (.9,1.5) -- (2.2,.2) -- (2,0) -- (1.4,-.6) -- (2.4,-1.6) -- (2.57,-1.43); 
	\draw[magenta, very thick] (2.63,-1.37) -- (3,-1) -- (2.83,-.83);
	\draw[magenta, very thick] (2.77,-.77) -- (2.7,-.7) -- (2.9,-.5) -- (3,-.6) -- (2.4,-1.2) -- (3,-1.8) -- (3.6,-1.2) -- (3.2,-.8) -- (3.65,-.35) -- (3.3,0) -- (2.9,.4) -- (3.97,1.47);
	\draw[magenta, very thick] (4.03,1.53) -- (4.25,1.75) -- (4,2) -- (3.75,1.75) -- (4.75,.75) -- (4.5,.5) -- (5,0) -- (5.7,-.7) -- (6.4,0) -- (7.1,.7) -- (7.37,0.43);
	\draw[magenta, very thick] (7.43,0.37) -- (8,-.2) -- (7.6,-.6) -- (7.53,-.53);
	\draw[magenta, very thick] (7.47,-.47) -- (7.2,-.2) -- (7.45,0.05) -- (7.75,-.25) -- (7.4,-.6) -- (6.9,-.1) -- (7.6,.6) -- (8.2,0) -- (8.6,-.4) -- (8.1,-.9) -- (8.5,-1.3) -- (9.3,-.5) -- (9,-.2) -- (9.2,0) -- (9.5,.3) -- (9.2,.6) -- (10.2,1.6) -- (10.5,1.3) -- (9.8,.6) -- (10.2,.2) -- (10,0);
	
	\node[circle, fill=black, inner sep=1.5pt, label=below:{\footnotesize$(0,0)$}] at (0.4,0) {};
	\node[circle, fill=black, inner sep=1.5pt, label=below:{\footnotesize$(n,0)$}] at (10,0) {};

	\node[magenta] at (.8,.67) {\footnotesize$-1$};
	\node[magenta] at (.8,.11) {\scriptsize$-2$};
	\node[magenta] at (2.65,-.28) {\footnotesize$+1$};
	\node[magenta] at (2.68,-1.1) {\scriptsize$+2$};
	\node[magenta] at (3.9,.67) {\footnotesize$-1$};
	\node[magenta] at (3.98,1.75) {\scriptsize$+1$};
	\node[magenta] at (5.68,-.28) {\footnotesize$+1$};
	
	\begin{pgfonlayer}{background}
	\fill[fill=magenta!20] (.4,0) -- (.8,.4) -- (1.2,0) -- (.3,-.9) -- (-.6,0) -- (.4,1) -- (.9,1.5) -- (2.2,.2) -- (2,0) -- (1.4,-.6) -- (2.4,-1.6) -- (3,-1) -- (2.7,-.7) -- (2.9,-.5) -- (3,-.6) -- (2.4,-1.2) -- (3,-1.8) -- (3.6,-1.2) -- (3.2,-.8) -- (3.65,-.35) -- (3.3,0) -- (2.9,.4) -- (4.25,1.75) -- (4,2) -- (3.75,1.75) -- (4.75,.75) -- (4.5,.5) -- (5,0) -- (5.7,-.7) -- (6.4,0) -- (7.1,.7) -- (8,-.2) -- (7.6,-.6) -- (7.2,-.2) -- (7.45,0.05) -- (7.75,-.25) -- (7.4,-.6) -- (6.9,-.1) -- (7.6,.6) -- (8.2,0) -- (8.6,-.4) -- (8.1,-.9) -- (8.5,-1.3) -- (9.3,-.5) -- (9,-.2) -- (9.2,0) -- (9.5,.3) -- (9.2,.6) -- (10.2,1.6) -- (10.5,1.3) -- (9.8,.6) -- (10.2,.2) -- (10,0);
	\fill[fill=magenta!40] (.4,0) -- (.8,.4) -- (1.2,0);
	\fill[fill=magenta!40] (2.6,-1.4) -- (3,-1) -- (2.8,-.8) -- (2.4,-1.2) -- (2.6,-1.4);
	\fill[fill=magenta!40] (7,0) -- (7.4,0.4) -- (7.8,0);
	\fill[fill=magenta!40] (7.5,-.5) -- (7.2,-.2) -- (7.45,0.05) -- (7.75,-.25) -- (7.5,-.5);
	\end{pgfonlayer}

	\draw[very thick, -latex] (4.3,-1.7) -- (2.5,-3);
	\draw[very thick, -latex] (5.7,-1.7) -- (7.5,-3);
	
	\begin{scope}[shift={(-1.3,-4)}, scale=.55]
		\draw[thick, lightgray] (.4,-1.6) -- (.4,1.7);
		\draw[thick, lightgray] (10,-1.6) -- (10,1.7);
		\draw[very thick, dashed, Purple] (.4,0) -- (.4,1);
		\draw[very thick, dashed, Purple] (10,0) -- (10,1.4);
		\draw[very thick, -latex] (-0.3,0) -- (11,0);
		
		\draw[Purple, very thick] (0.4,1) -- (.9,1.5) -- (2.2,.2) -- (2,0) -- (1.4,-.6) -- (2.4,-1.6) -- (2.6,-1.4) -- (3,-1.8) -- (3.6,-1.2) -- (3.2,-.8) -- (3.65,-.35) -- (3.3,0) -- (2.9,.4) -- (4,1.5) -- (4.75,.75) -- (4.5,.5) -- (5,0) -- (5.7,-.7) -- (6.4,0) -- (7.1,.7) -- (7.4,.4) -- (7.6,.6) -- (8.2,0) -- (8.6,-.4) -- (8.1,-.9) -- (8.5,-1.3) -- (9.3,-.5) -- (9,-.2) -- (9.2,0) -- (9.5,.3) -- (9.2,.6) -- (10,1.4);
		
		\node[circle, fill=black, inner sep=1.2pt] at (0.4,0) {};
		\node[circle, fill=black, inner sep=1.2pt] at (10,0) {};
		
		\node[Purple] at (1.05,.5) {\footnotesize$-1$};
		\node[Purple] at (2.65,-.4) {\footnotesize$+1$};
		\node[Purple] at (3.9,.5) {\footnotesize$-1$};
		\node[Purple] at (5.68,-.28) {\scriptsize$+1$};
		
		\begin{pgfonlayer}{background}
			\fill[fill=Purple!20] (0.4,0) -- (0.4,1) -- (.9,1.5) -- (2.2,.2) -- (2,0) -- (1.4,-.6) -- (2.4,-1.6) -- (2.6,-1.4) -- (3,-1.8) -- (3.6,-1.2) -- (3.2,-.8) -- (3.65,-.35) -- (3.3,0) -- (2.9,.4) -- (4,1.5) -- (4.75,.75) -- (4.5,.5) -- (5,0) -- (5.7,-.7) -- (6.4,0) -- (7.1,.7) -- (7.4,.4) -- (7.6,.6) -- (8.2,0) -- (8.6,-.4) -- (8.1,-.9) -- (8.5,-1.3) -- (9.3,-.5) -- (9,-.2) -- (9.2,0) -- (9.5,.3) -- (9.2,.6) -- (10,1.4) -- (10,0);
		\end{pgfonlayer}
	\end{scope}
	
	\begin{scope}[shift={(5.7,-4)}, scale=.55]
		\draw[thick, lightgray] (.4,-1.6) -- (.4,1.7);
		\draw[thick, lightgray] (10,-1.6) -- (10,1.7);
		\draw[very thick, dashed, orange] (.4,0) -- (.4,1);
		\draw[very thick, dashed, orange] (10,0) -- (10,1.4);
		\draw[very thick, -latex] (-1.1,0) -- (11,0);
		
		\draw[dotted, thin] (0.4,1) -- (.9,1.5) -- (2.2,.2) -- (2,0) -- (1.4,-.6) -- (2.4,-1.6) -- (2.6,-1.4) -- (3,-1.8) -- (3.6,-1.2) -- (3.2,-.8) -- (3.65,-.35) -- (3.3,0) -- (2.9,.4) -- (4,1.5) -- (4.75,.75) -- (4.5,.5) -- (5,0) -- (5.7,-.7) -- (6.4,0) -- (7.1,.7) -- (7.4,.4) -- (7.6,.6) -- (8.2,0) -- (8.6,-.4) -- (8.1,-.9) -- (8.5,-1.3) -- (9.3,-.5) -- (9,-.2) -- (9.2,0) -- (9.5,.3) -- (9.2,.6) -- (10,1.4);
		
		\draw[orange, very thick] (.4,0) -- (.8,.4) -- (1.2,0) -- (.3,-.9) -- (-.6,0) -- (.4,1);
		\draw[orange, very thick] (2.6,-1.4) -- (3,-1) -- (2.85,-.85);
		\draw[orange, very thick] (2.75,-.75) -- (2.7,-.7) -- (2.9,-.5) -- (3,-.6) -- (2.4,-1.2) -- (2.6,-1.4);
		\draw[orange, very thick] (4,1.5) -- (4.25,1.75) -- (4,2) -- (3.75,1.75) -- (4,1.5);
		\draw[orange, very thick] (7.4,0.4) -- (8,-.2) -- (7.6,-.6) -- (7.55,-.55);
		\draw[orange, very thick] (7.45,-.45) -- (7.2,-.2) -- (7.45,0.05) -- (7.75,-.25) -- (7.4,-.6) -- (6.9,-.1) -- (7.4,.4);
		\draw[orange, very thick] (10,1.4) -- (10.2,1.6) -- (10.5,1.3) -- (9.8,.6) -- (10.2,.2) -- (10,0);
		
		\node[circle, fill=black, inner sep=1.2pt] at (0.4,0) {};
		\node[circle, fill=black, inner sep=1.2pt] at (10,0) {};
		
		\node[orange] at (0.3,-.4) {\scriptsize$-1$};
		
		\begin{pgfonlayer}{background}
			\fill[fill=orange!20] (.4,0) -- (.8,.4) -- (1.2,0) -- (.3,-.9) -- (-.6,0) -- (.4,1) -- cycle;
			\fill[fill=orange!20] (2.6,-1.4) -- (3,-1) -- (2.7,-.7) -- (2.9,-.5) -- (3,-.6) -- (2.4,-1.2) -- (2.6,-1.4);
			\fill[fill=orange!20] (4,1.5) -- (4.25,1.75) -- (4,2) -- (3.75,1.75) -- (4,1.5);
			\fill[fill=orange!20] (7.4,0.4) -- (8,-.2) -- (7.6,-.6) -- (7.2,-.2) -- (7.45,0.05) -- (7.75,-.25) -- (7.4,-.6) -- (6.9,-.1) -- (7.4,.4);
			\fill[fill=orange!40] (7.5,-.5) -- (7.2,-.2) -- (7.45,0.05) -- (7.75,-.25) -- (7.5,-.5);
			\fill[fill=orange!20] (10,1.4) -- (10.2,1.6) -- (10.5,1.3) -- (9.8,.6) -- (10.2,.2) -- (10,0);
		\end{pgfonlayer}
	\end{scope}
\end{tikzpicture} \vspace*{-4mm}
\captionsetup{margin=5mm, font=footnotesize}

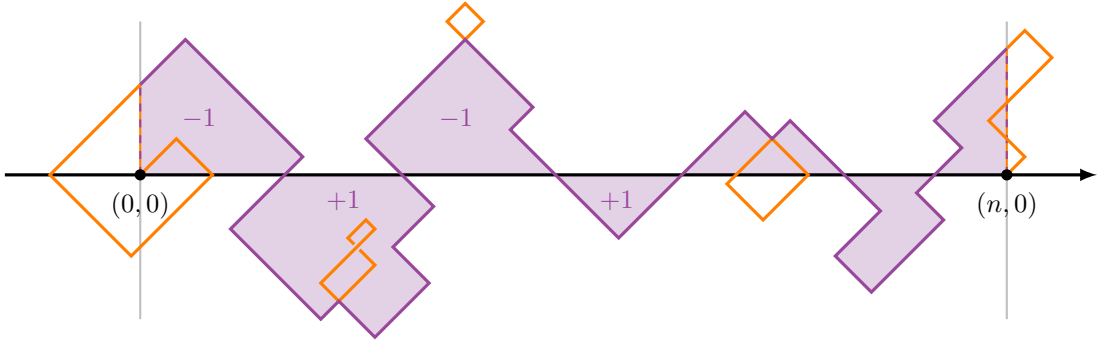
\captionof{figure}{The decomposition of a word $w$. The purple curve is obtained from the interval between the last crossing of the line $x=0$ and the first crossing of $x=n$, after loop-erasure.}
\end{center}
By \say{contribution}, we mean that the winding number distribution of $w$ splits as the sum of the winding numbers of the purple curve and the orange curves (with the added dashed segments). Integrating against the $y$-coordinate gives two contributions to $B_y(\overline w)$, which we denote $B_y(\text{purple curve})$ and $B_y(\text{orange curves})$.

\bigbreak

$\blacktriangleright$ We first estimate the contribution of the purple main path. Note that this path is simple (and we cured out boundary mess), so the winding number of any point is $\pm 1$ if the point lies between the $x$-axis and the curve - more precisely $+1$ if it lies below the $x$-axis and $-1$ if it lies above - and $0$ otherwise. In particular, the contribution is non-positive, and bounded by the contribution of the following green area:
\begin{center}
		\begin{tikzpicture}
	\clip (-.5,-1.9) rectangle (15,1.9);
	\draw[TealBlue!30, fill=TealBlue!30] (0.4,1) -- (.9,1.5) -- (.9,0) -- (.4,0)
	(2,-1.2) -- (2.4,-1.6) -- (2.4,0) -- (2,0)
	(2.6,-1.4) -- (3,-1.8) -- (3,0) -- (2.6,0)
	(3.2,-.8) -- (3.3,-.7) -- (3.3,0) -- (3.2,0)
	(4,1.5) -- (4.5,1) -- (4.5,0) -- (4,0)
	(5,0) -- (5.7,-.7) -- (5.7,0)
	(6.4,0) -- (7.1,.7) -- (7.1,0)
	(7.4,.4) -- (7.6,.6) -- (7.6,0) -- (7.4,0)
	(8.2,0) -- (8.6,-.4) -- (8.6,0)
	(9,-.2) -- (9.2,0) -- (9,0)
	(9.5,.9) -- (10,1.4) -- (10,0) -- (9.5,0);
	\draw[JungleGreen!40, fill=TealBlue!50] (.9,1.5) -- (2,.4) -- (2,0) -- (.9,0)
	(2.4,-1.6) -- (2.6,-1.4) -- (2.6,0) -- (2.4,0)
	(3,-1.8) -- (3.2,-1.6) -- (3.2,0) -- (3,0)
	(3.3,.8) -- (4,1.5) -- (4,0) -- (3.3,0)
	(4.5,.5) -- (5,0) -- (4.5,0)
	(5.7,-.7) -- (6.4,0) -- (5.7,0)
	(7.1,.7) -- (7.4,.4) -- (7.4,0) -- (7.1,0)
	(7.6,.6) -- (8.2,0) -- (7.6,0)
	(8.6,-1.2) -- (9,-.8) -- (9,0) -- (8.6,0)
	(9.2,0) -- (9.5,.3) -- (9.5,0);
	
	\draw[Purple, dashed] (0.4,1) -- (.9,1.5) -- (2.2,.2) -- (2,0) -- (1.4,-.6) -- (2.4,-1.6) -- (2.6,-1.4) -- (3,-1.8) -- (3.6,-1.2) -- (3.2,-.8) -- (3.65,-.35) -- (3.3,0) -- (2.9,.4) -- (4,1.5) -- (4.75,.75) -- (4.5,.5) -- (5,0) -- (5.7,-.7) -- (6.4,0) -- (7.1,.7) -- (7.4,.4) -- (7.6,.6) -- (8.2,0) -- (8.6,-.4) -- (8.1,-.9) -- (8.5,-1.3) -- (9.3,-.5) -- (9,-.2) -- (9.2,0) -- (9.5,.3) -- (9.2,.6) -- (10,1.4);

	\draw[Purple, very thick] (0.4,1) -- (.9,1.5) -- (2,.4)
	(2,-1.2) -- (2.4,-1.6) -- (2.6,-1.4) -- (3,-1.8) -- (3.2,-1.6)
	(3.2,-.8) -- (3.3,-.7)
	(3.3,.8) -- (4,1.5) -- (4.5,1)
	(4.5,.5) -- (5,0) -- (5.7,-.7) -- (6.4,0) -- (7.1,.7) -- (7.4,.4) -- (7.6,.6) -- (8.2,0) -- (8.6,-.4)
	(8.6,-1.2) -- (9,-.8)
	(9,-.2) -- (9.2,0) -- (9.5,.3)
	(9.5,.9) -- (10,1.4);

	\begin{scope}[shift={(6.5,-1.1)}, scale=1.7]
		\clip (3.9,.65) circle (10mm);
	
		\draw[Purple, dashed] (3.2,-.8) -- (3.65,-.35) -- (3.3,0) -- (2.9,.4) -- (4,1.5) -- (4.75,.75) -- (4.5,.5);
		\draw[TealBlue!30, fill=TealBlue!30] (4,1.5) -- (4.5,1) -- (4.5,0) -- (4,0)
			(2.8,-.8) -- (3.3,-.7) -- (3.3,0) -- (2.8,0);
		\draw[TealBlue!50, fill=TealBlue!50] (3.3,.8) -- (4,1.5) -- (4,0) -- (3.3,0)
			(4.5,.5) -- (5,0) -- (4.5,0);
	
		\draw[Purple, very thick] (3.3,.8) -- (4,1.5) -- (4.5,1)
			(4.5,.5) -- (5,0);
	
		\fill[NavyBlue!40] (3.3,0) -- (4,.7) -- (4,0);
		\fill[NavyBlue!20] (4,0) -- (4,.5) -- (4.5,0);
		\draw[NavyBlue, thick] (3.3,0) -- (4,.7);
		\draw[NavyBlue, thick] (4,.5) -- (4.5,0);
	
		\draw[very thick] (2.8,0) -- (5,0);
		
		\draw[black!80, decorate,decoration={brace,amplitude=3pt},yshift=-1pt] (3.98,0) -- (3.32,0) 
		node [black,midway,yshift=-3mm] {\color{black!80}\footnotesize $a_i$};
		\node[circle, fill=NavyBlue, inner sep=.8pt] at (3.76,.23) {};
		\node[circle, fill=NavyBlue, inner sep=.8pt] at (4.16,.16) {};
	\end{scope}
	
	\draw[lightgray, very thin] (3.9,1.65) -- (13.1,1.705);
	\draw[lightgray, very thin] (3.75,-.34) -- (13,-1.69);
	\draw[gray, thin] (3.9,.65) circle (10mm);
	\draw[gray, thin, shift={(6.5,-1.1)}, scale=1.7] (3.9,.65) circle (10mm);
	
	\draw[thick, gray!80] (.4,-1.6) -- (.4,1.7);
	\draw[thick, gray!80] (10,-1.6) -- (10,1.7);
	\draw[very thick, -latex] (-.4,0) -- (11,0);
	
	\node[circle, fill=black, inner sep=1.5pt] at (0.4,0) {};
	\node[circle, fill=white, fill opacity=.8, inner sep=.5pt] at (0.4,-.3) {\footnotesize$0$};
	\node[circle, fill=black, inner sep=1.5pt] at (10,0) {};
	\node[circle, fill=white, fill opacity=.8, inner sep=.5pt] at (10,-.3) {\footnotesize$n$};
\end{tikzpicture}
		\captionsetup{margin=5mm, font=footnotesize}
		
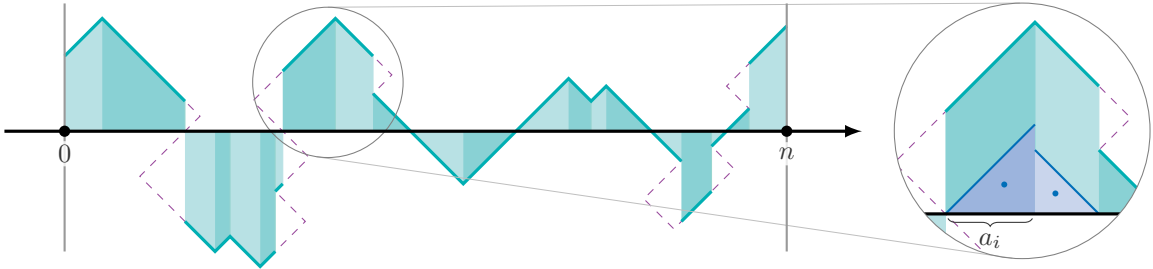
\captionof{figure}{The green area decomposed into $k'$ trapezoids, and the blue triangles included inside.}
\end{center}
The green area is composed of $k'$ trapezoids/triangles, with $k'\le 2k(w)-2$ (each segment delimits at most $2$ trapezoids if it crosses the $x$-axis, at most $1$ otherwise, and the first and last segments cannot cross the axis). In turn, we bound the contribution of each slice by that of a triangle (with basis and height $a_i$) included inside it:
\begin{equation} \label{eq:purple_contribution}
 -B_y(\text{purple curve}) \ge -B_y(\text{green zone}) \ge -B_y(\text{blue triangles}) = \sum_{i=1}^{k'} \frac{a_i^2}2\cdot\frac{a_i}3.
\end{equation}
Finally, since $\sum_{i=1}^{k'}a_i=n$, the generalized mean inequality gives
\[ -B_y(\text{purple curve})\ge \frac16\frac{n^3}{k'^2} \ge \frac1{24}\frac{n^3}{(k(w)-1)^2}.\]

$\blacktriangleright$ To control the contribution of orange curves, we need to control both their total area and the $y$-coordinates. Observe that the purple curve (without the dashed segments) has length at least $n$ as it joins the lines $x=0$ and $x=n$. It follows that the orange curves have total length at most $\Delta$, hence the two dashed segments too. Therefore
\begin{align} 
	B_y(\text{orange curves})
	& \le \iint_{\RR^2} \abs{y}\cdot \abs{W_{\text{orange curves}}(x,y)} \cdot \dif x \dif y \notag \\
	& \le y_{\max} \cdot \iint_{\RR^2}\abs{W_{\text{orange curves}}(x,y)} \cdot \dif x \dif y \notag \\
	& \le y_{\max}\cdot I\cdot (2\Delta)^2, \label{eq:orange_contribution}
\end{align}
where $y_{\max}$ is the largest distance from points of the path $w$ to the $x$-axis, and $I$ is the isoperimetric constant of $(\RR^2,\norm{\,\cdot\,}_{\Mink,X})$ (here $I=\frac18$). If $y_{\max}$ is reasonably small, say $y_{\max}\le L\max\big\{\Delta;\, \frac n{k(w)-1}\big\}$ for $L=10\max\{\sqrt I,1\}$, then
\begin{align*}
-B_y(\bar w)  \ge -B_y(\text{purple}) -B_y(\text{orange}) \ge  \frac1{24}\frac{n^3}{(k(w)-1)^2} - 4LI\cdot \Delta^2 \cdot \max\left\{\Delta;\,\frac{n}{k(w)-1}\right\}.
\end{align*}
\bigbreak
$\blacktriangleright$ The only remaining case is when $y_{\max}$ is unreasonably large: $y_{\max}\ge L\max\big\{\Delta;\, \frac n{k(w)-1}\big\}$. In particular, $y_{\max}-\Delta$ is larger than $(1-\frac1L)y_{\max}$, $(L-1)\Delta$ and $(L-1)\frac n{k(w)-1}$.

\newpage

We improve our bound on $-B_y(\text{green zone})$ using that the curve $w$ goes through some point $p=(x_p,\pm y_{\max})$ far away from the $x$-axis: 
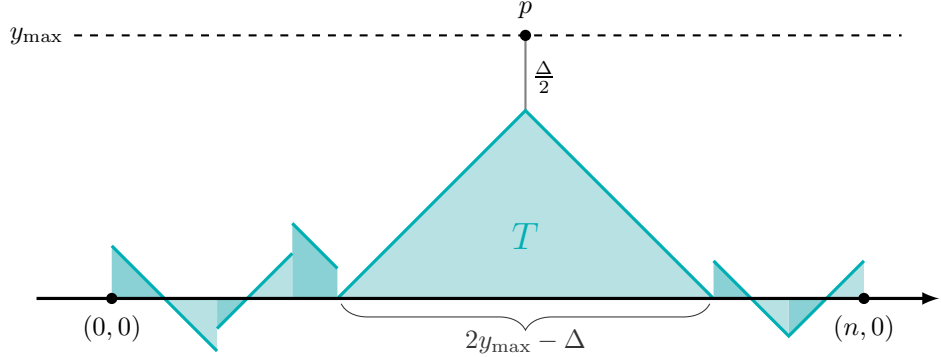
\vspace*{-2mm}
\begin{center}
\begin{tikzpicture}[scale=.9]
	\draw[line width=.8pt, dashed] (-.5,3.5) -- (10.5,3.5);
	\node at (-1,3.5) {\footnotesize$y_{\max}$};
	\draw[thick, gray] (5.5,3.5) -- (5.5,2.5);
	\node at (5.75,2.95) {\footnotesize$\Delta$};
	\draw[black!80, decorate,decoration={brace,amplitude=10pt},yshift=-2pt] (7.95,0) -- (3.05,0) 
		node [black,midway,yshift=-5mm] {\color{black!80}\footnotesize $2(y_{\max}-\Delta)$};
	\node[circle, fill=black, inner sep=1.5pt, label=above:{\footnotesize$p$}] at (5.5,3.5) {};
	\draw[very thick, TealBlue, fill=TealBlue!30] (3,0) -- (5.5,2.5) --  (8,0);
	\node[TealBlue] at (5.5,.83) {\large$T$};
	
	\fill[fill=TealBlue!50] (0,.7) -- (.7,0) -- (0,0)
	(1.4,-.4) -- (1.8,0) -- (1.4,0)
	(2.4,1) -- (3,.4) -- (3,0) -- (2.4,0)
	(8,.5) -- (8.5,0) -- (8,0)
	(9,-.5) -- (9.5,0) -- (9,0);
	
	\fill[fill=TealBlue!30] (.7,0) -- (1.4,-.7) -- (1.4,0) -- (.7,0)
	(1.8,0) -- (2.4,.6) -- (2.4,0)
	(8.5,0) -- (9,-.5) -- (9,0)
	(9.5,0) -- (10,.5) -- (10,0);
	\draw[very thick, Purple] (0,.7) -- (.7,0) -- (1.4,-.7)
	(1.4,-.4) -- (1.8,0) -- (2.4,.6)
	(2.4,1) -- (3,.4)
	(8,.5) -- (8.5,0) -- (9,-.5) -- (9.5,0) -- (10,.5);
	
	\draw[very thick, -latex] (-1,0) -- (11,0);
	\node[circle, fill=black, inner sep=1.5pt, label=below:{\footnotesize$(0,0)$}] at (0,0) {};
	\node[circle, fill=black, inner sep=1.5pt, label=below:{\footnotesize$(n,0)$}] at (10,0) {};
\end{tikzpicture} \vspace*{-3mm}
\end{center}
In this case, there exists a large triangle $T$ which the $\RR$-word $w$ cannot enter. Indeed, for any point $q$ in the interior of $T$, we have
\begin{align*}
d_\Mink\big((0,0),q\big)+d_\Mink(q,p) & >d_\Mink\big((0,0),p\big)+\Delta \\
d_\Mink(p,q)+d_\Mink\big(q,(n,0)\big) & >d_\Mink\big(p,(n,0)\big)+\Delta
\end{align*}
so any $\RR$-word passing through the four points $(0,0)$, $p$, $q$ and $\hat w=(n,0)$ would have length $>n+\Delta$. This implies that the triangle $T$ must be included inside the green region. Moreover, as previously, we have a lower bound on $-B_y$ for the green regions on both sides of $T$, composed of at most $2k(w)-2$ trapezoids.

Combining equations (\ref{eq:purple_contribution}), (\ref{eq:orange_contribution}) and the existence of $T$, we get
\begin{align*}
-B_y(\bar w)
& \ge -B_y(\text{green area}) - y_{\max} \cdot 4I\Delta^2 \\
& \ge \frac1{24}\frac{\big(n-2(y_{\max}-\Delta)\big)^3}{(k(w)-1)^2} + \frac13(y_{\max}-\Delta)^3 -y_{\max}\cdot 4I\Delta^2 \\
& \ge \frac1{24}\frac{n^3}{(k-1)^2} - \frac14(y-\Delta)\frac{n^2}{(k-1)^2}+\frac13(y-\Delta)^3 -4I\cdot y\Delta^2 \\
& \ge \frac1{24}\frac{n^3}{(k-1)^2} + (y-\Delta)^3 \left(\frac13- \frac1{4(L-1)^2} -\frac{4IL}{(L-1)^3} \right) \\
& \ge \frac1{24}\frac{n^3}{(k(w)-1)^2}.
\end{align*}
(In the third step, we use $(n-h)^3\ge n^3-3n^2h$ for all $h\le 3n$ and $y_{\max}\le \frac12(n+\Delta)$.)
\end{proof}

We deduce the following quantified version of the main result of \cite{Stoll_3step}:
\begin{cor}\label{rem:BrLD}
	There exists $C'>0$ such that, for all $w\in X_\RR^*$ representing $g_n$, we have
	\[ \ell(w)-\norm{g_n}_{\Stoll,X}\ge C'\cdot k(w)^{-\frac23}\cdot \norm{g_n}_{\Stoll,X}. \vspace*{-1mm}\]
	If moreover $n\equiv 0\pmod 6$, then $g_n\in \Ec$ and $\norm{g_n}_{X} - \norm{g_n}_{\Stoll,X}= \Omega\left( \norm{g_n}_{\Stoll,X}^{\frac13}\right)$.
\end{cor}
\begin{proof}
Observe that $\hat w=\hat g_n=(n,0)$ and $B_y(\bar w)=B_y(g_n)=0$, therefore
	\begin{align*}
	 & 0  \ge \frac1{24}\frac{n^3}{(k(w)-1)^2} - C\Delta^2\max\left\{\Delta;\,\frac n{k(w)-1}\right\} \\
\iff & \max\left\{ \Delta^3;\, \frac{n\Delta^2}{k(w)-1}\right\} \ge \frac1{24C}\frac{n^3}{(k(w)-1)^2} \\
\iff & \Delta \ge  \min\left\{ \sqrt[3]{\frac1{24C}\frac{n^3}{(k(w)-1)^2}};\, \sqrt{\frac1{24C}\frac{n^2}{k(w)-1}} \right\}.
	\end{align*}
	Since $\norm{g_n}_{\Stoll,X}=n$, this implies that \vspace*{-3mm}
	\[ \ell(w)-\norm{g_n}_{\Stoll,X} \ge \sqrt{\frac1{24C}} \cdot k(w)^{-\frac23} \cdot \norm{g_n}_{\Stoll,X}.\vspace*{-1mm} \]
	Finally, as $k(w)\le \ell(w)$ for all $w\in X^*$, we have $\norm{g_n}_X -\norm{g_n}_{\Stoll,X} \ge C' \cdot \norm{g_n}_{\Stoll,X}^\frac13$.
\end{proof}

\medbreak

The first part of Corollary \ref*{rem:BrLD} matches the best known upper bound: \vspace*{-7.5mm}
\begin{adjustwidth}{2mm}{2mm}
\begin{lemma*}[{\cite[Lemma 40]{Gianella}}]
	Let $\bar H$ be a simply connected $s$-step nilpotent Lie group, and $X$ a finite Lie generating set. There exists $C''>0$ such that, for every $K\gg 1$ and every $g\in\bar H$, there exists an $\RR$-word $w\in X_\RR^*$ representing $g$ such that \vspace*{-1mm}
	\[ k(w)\le K \quad\text{and}\quad \ell(w)-\norm{g}_{\Stoll,X} \le C'' \cdot K^{-\frac2s} \cdot \norm{g}_{\Stoll,X} +C''. \vspace*{-2mm}\]
\end{lemma*}
\end{adjustwidth}
Whether this bound is sharp in $s$-step nilpotent groups with $s\ge 4$ remains open. The second part of Corollary \ref*{rem:BrLD} disproves a conjecture of Breuillard and Le Donne \cite[Conjecture 6.5]{BreuillardLeDonne}. (The conjecture states that the difference should be an $O(1)$.)

\subsection{Matching lower bound in a virtually Engel group}

In this paragraph, we prove Theorem \ref{thm:inter}, considering the semi-direct product
\[ \vEc = \Ec \rtimes C_2 = \la a,t \;\big|\; t^2=1;\; [a,[a,a^t]]=[a^t,[a,a^t]]\text{ commutes with }a,a^t \ra\]
(so $C_2=\la t\ra$ acts by symmetry along the $y$-axis, and in particular $tat=b^{-1}$), with the generating set $S=\{a^\pm,t\}$. First, we may compute $A(S)$ and $P(S)$, so that Theorem 1 gives the upper bound $\gamma_\geod(n)\preceq \exp\!\big(n^{3/5}\cdot \log(n)\big)$.
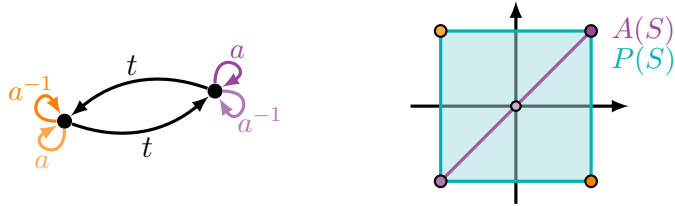
\begin{center}
	\begin{tikzpicture}[scale=1]
	\node[circle, fill=black, inner sep=2pt] (0) at (1,.2) {};
	\node[circle, fill=black, inner sep=2pt] (1) at (-1,-.2) {};
	\draw[very thick, -latex, in=40, out=160] (0) to (1) node [xshift=-1mm, yshift=1mm, pos=.45, label=$t$] {};
	\draw[very thick, -latex, in=-140, out=-20] (1) to (0) node [xshift=1mm, yshift=-1mm, pos=.45, label=below:$t$] {};
	
	\draw[very thick, Purple, -latex] (0) to[out=90, in=30, loop] (0) {};
	\draw[very thick, Purple!70, -latex] (0) to[out=0, in=-60, loop] (0) {};
	\draw[very thick, orange!70, -latex] (1) to[out=-90, in=-150, loop] (1) {};
	\draw[very thick, orange, -latex] (1) to[out=180, in=120, loop] (1) {};
	\node[Purple] at (1.3,.72) {\small$a$};
	\node[Purple!70] at (1.6,-.18) {\small$a^{-1}$};
	\node[orange!70] at (-1.3,-.74) {\small$a$};
	\node[orange] at (-1.45,.24) {\small$a^{-1}$};
	
	\begin{scope}[shift={(5,0)}]
	\draw[very thick, -latex] (-1.4,0) -- (1.5,0);
	\draw[very thick, -latex] (0,-1.3) -- (0,1.4);
	\draw[very thick, TealBlue, fill=TealBlue!40, fill opacity=.5] (1,1) -- (-1,1) -- (-1,-1) -- (1,-1) -- cycle;
	\draw[very thick, Purple!90] (1,1) -- (-1,-1);
	\draw[thick, black, fill=Purple] (1,1) circle (2.2pt);
	\draw[thick, black, fill=Purple!40] (0,0) circle (1.8pt);
	\draw[thick, black, fill=Purple!70] (-1,-1) circle (2.2pt);
	\draw[thick, black, fill=orange!70] (-1,1) circle (2.2pt);
	\draw[thick, black, fill=orange] (1,-1) circle (2.2pt);
	\node[Purple] at (1.7,1) {\small$A(S)$};
	\node[TealBlue] at (1.7,.6) {\small$P(S)$};
	\end{scope}
	\end{tikzpicture}
\end{center}
It remains to prove a matching lower bound, i.e., to construct a lot of geodesics. We fix $\kappa>0$ a small constant (to be determined) and $0<\varepsilon<\frac1{10}$. For any $n$ even integer, let $K\approx \kappa n^{3/5}$ be another even integer, and $m=\frac n{2K}$. We show that words of the form
\[ w = a^{m_1}ta^{-(m_1+m_2)}ta^{m_2+m_3}t\ldots ta^{-(m_{K-1}+m_K)}ta^{m_K} \in S^*\]
with $2\sum_i m_i=n$ and $\abs{m_i-m}\le n^\varepsilon$, are all geodesics for $n$ large enough.
\begin{center}
	\begin{tikzpicture}[scale=.4,
	every loop/.style={},
	lefthalf/.style={semicircle, anchor=chord center, rotate=90, minimum size=2pt, inner sep=0},		righthalf/.style={semicircle, anchor=chord center, rotate=-90, minimum size=2pt, inner sep=0}]
		
	\newcommand{\PtoO}[2]{
		\node[lefthalf, draw=Purple, fill=Purple] at (#1,#2) {};
		\node[righthalf, draw=orange, fill=Orange] at (#1,#2) {};
		\node[minimum size=3pt, inner sep=0] at (#1,#2) {} edge[very thick, black, out=45, in=135, loop] ();}
		
	\newcommand{\OtoP}[2]{
		\node[lefthalf, draw=orange, fill=orange] at (#1,#2) {};
		\node[righthalf, draw=Purple, fill=Purple] at (#1,#2) {};
		\node[minimum size=3pt, inner sep=0] at (#1,#2) {} edge[very thick, black, out=-45, in=-135, loop] ();}
		
	\begin{scope}
		\clip (-.7,-4.7) rectangle (30.7,5.7);
		\draw[very thin, lightgray, cm={1,1,1,-1,(0,0)}] (-4,-4) grid (20,20);
	\end{scope}
	
	\draw[very thick, -latex] (-1,0) -- (31.2,0);
	\draw[very thick, -latex] (0,-5) -- (0,6);

	\draw[ultra thick, Purple] (0,0) -- (2,2);
	\draw[ultra thick, orange] (2,2) -- (6,-2);
	\draw[ultra thick, Purple] (6,-2) -- (12,4);
	\draw[ultra thick, orange] (12,4) -- (18,-2);
	\draw[ultra thick, Purple] (18,-2) -- (22,2);
	\draw[ultra thick, orange] (22,2) -- (27,-3);
	\draw[ultra thick, Purple] (27,-3) -- (30,0);
	
	\draw[thick, black] (12,0) -- (12,4);
	\node at (12.8,1.6) {\footnotesize$m_3$};
	\node at (12.5,5) {\footnotesize$t$};
	\node at (7.8,-1.3) {\footnotesize$a^{m_2}$};
	\node at (9.3,2.2) {\footnotesize$a^{m_3}$};
	
	\node[circle, fill=Purple, minimum size=4pt, inner sep=0] at (0,0) {};
	\node[circle, fill=Purple, minimum size=4pt, inner sep=0] at (30,0) {};
	
	\PtoO{2}{2}
	\OtoP{6}{-2}
	\PtoO{12}{4}
	\OtoP{18}{-2}
	\PtoO{22}{2}
	\OtoP{27}{-3}
\end{tikzpicture}	 \vspace*{-2mm}
	\captionsetup{margin=5mm, font=footnotesize}
	\captionof{figure}{The path in $\Ec$ corresponding to a geodesic in $\vEc$.}
\end{center}
$\blacktriangleright$ First, we compute $B_y$ for the corresponding element $\overline w$. Note that $\overline w\in\Ec$ (as $K$ is even), so this makes sense. Let $\delta_i=m_i-m$ (so that $\sum_i\delta_i=0$), we have
\begin{equation} \label{eq:B_of_w}
	\newcommand{\se}{\hspace*{-.25mm}=\hspace*{-.25mm}} 
	-B_y(\overline w)
	\se \sum_{i=1}^K \frac{m_i^3}3
	\se \sum_{i=1}^K \frac{(m+\delta_i)^3}3
	\se \frac{Km^3}3 + \sum_{i=1}^K \!\Big(\hspace*{-.3mm}m\delta_i^2 + \frac{\delta_i^3}3\Big)
	\se \frac{n^3}{24K^2} + O(n^{1+2\varepsilon}).
\end{equation}

$\blacktriangleright$ Next, we take a shorter word $v\in\{a^\pm,t\}^*$ ending up in the same coset as $w$ (that is the $\Ec$ coset), and with same endpoint $\hat v=\hat w=(n,0)$, and prove that $B_y(\bar v)<B_y(\bar w)$. This implies that $w$ is a geodesic, as no shorter word represents the same element.

Consider the decomposition $\tilde v\in X^*$. Its coarse length and length are
\begin{align*}
k(\tilde v)-1 & \le \text{\say{number of $t$ in $v$}} \hspace{6.5mm} = K-d\quad \text{(for some even $d$)} \\
\ell(\tilde v) & = \text{\say{number of $a^{\pm1}$ in $v$}} \hspace{2mm}< n+d.
\end{align*}
Note that, as you need at least $n$ letters \say{$a^{\pm 1}$} to reach $\hat v=(n,0)$, the shortening relative to $w$ has to come from the number of \say{$t$} in $v$, so that $2\le d\le K-2$.
\medbreak
Combining Proposition \ref{prop:back_from_depths} and the equation (\ref{eq:B_of_w}), we get
\[ B_y(\bar w) - B_y(\bar v)
 \ge \frac{n^3}{24}\left(\frac1{(K-d)^2}-\frac1{K^2}\right) - C\cdot d^2 \cdot \max\left\{d;\frac n{K-d}\right\} - O(n^{1+2\varepsilon}). \]

Now we split into two cases:
\begin{itemize}[leftmargin=8mm]
	\item If $d\ge K-\frac nK$, then the term $\frac{n^3}{24(K-d)^2}$ dominates and $B_y(\bar w)>B_y(\bar v)$.
	
	
	\item If $d\le K-\frac nK$. We have $\frac1{(K-d)^2}-\frac1{K^2}\ge \frac {2d}{K^3}$ by the mean value theorem, hence \vspace*{-1mm}
	\begin{align*}
	B_y(\bar w) - B_y(\bar v)
	& \ge \frac{n^3}{12}\cdot \frac{d}{K^3} - C \cdot d\cdot K^2-O(n^{1+2\varepsilon}) \\ 
	& = d\left(\frac1{12\kappa^3}-C\kappa^2\right)\cdot n^{6/5} - o(n^{6/5}) \\[-1mm]
	& > 0 
	\end{align*}\vspace*{-9mm}
	
	as long as $\kappa<\sqrt[5]{\frac1{12C}}$ and $n$ is large enough.
\end{itemize}

\medbreak

$\blacktriangleright$ In conclusion, all those words $w$ are geodesics. Moreover there are many degrees of freedom left, leading to plenty of geodesics of length $n+K$. Indeed, the values of the partial sums $\big(\!\sum_{i=1}^r m_i\big)_{r=1,\ldots,K-1}$ can be picked independently in $[rm-\frac12 n^\varepsilon,rm+\frac12n^\varepsilon]$. This shows that, for all $n$ large enough even and $K\approx\kappa n^{3/5}$ even, we have
\[ \gamma_\geod\big(n+K\big) \ge (n^\varepsilon)^{K-1} \asymp \exp\!\big(\varepsilon\kappa \cdot n^{3/5}\cdot\log(n)\big). \]
Since the function $\gamma_\geod(n)$ is increasing, this bound extends for all $n$. \hfill$\square$

\section{Further questions} \label{sec:Qs}

Plenty of questions remain open. The question we would most like to see solved is
\begin{adjustwidth}{4mm}{4mm}
\textbf{Conjecture A.} If the geodesic growth of $(G,S)$ is polynomial, with $S$ a \emph{symmetric} generating set, then $G$ is virtually $2$-step nilpotent.
\end{adjustwidth}
Among the possible counter-examples (all virtually nilpotent), treating the virtually $3$-step nilpotent cases would be sufficient, as $G$ factors onto $G/\gamma_4(H)$. We emphasize \say{symmetric} as we have the following intriguing example:
\begin{adjustwidth}{4mm}{4mm}
	\textbf{Conjecture B.} The geodesic growth of $\vEc$ w.r.t.\ $S=\{a,b,(ab)^{-1},t\}$ is polynomial.
\end{adjustwidth}
In contrast, in an $\Ec$-by-finite group, any \emph{symmetric} generating set $S$ such that $P(S)$ has vertices on the $x$-axis will yield exponential geodesic growth (as the $x$-axis is fixed by automorphisms of $\Ec$), and any generating set $S$ such that $P(S)$ has no vertex on the $x$-axis should yield super-polynomial geodesic growth.
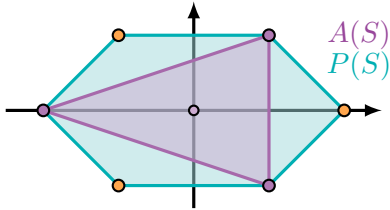
\vfill
\begin{center}
\begin{tikzpicture}[scale=1.2]
	\draw[very thick, -latex] (-2.5,0) -- (2.5,0);
	\draw[very thick, -latex] (0,-1.3) -- (0,1.45);
	\draw[very thick, TealBlue, fill=TealBlue!40, fill opacity=.5] (1,1) -- (-1,1) -- (-2,0) -- (-1,-1) -- (1,-1) -- (2,0) -- cycle;
	\draw[very thick, Purple!80, fill=Purple!40, fill opacity=.5] (1,1) -- (1,-1) -- (-2,0) -- cycle;
	\draw[thick, black, fill=Purple!70] (1,1) circle (2.2pt);
	\draw[thick, black, fill=Purple!70] (1,-1) circle (2.2pt);
	\draw[thick, black, fill=Purple!70] (-2,0) circle (2.2pt);
	\draw[thick, black, fill=Purple!40] (0,0) circle (1.8pt);
	
	\draw[thick, black, fill=TealBlue!70] (-1,1) circle (2.2pt);
	\draw[thick, black, fill=TealBlue!70] (-1,-1) circle (2.2pt);
	\draw[thick, black, fill=TealBlue!70] (2,0) circle (2.2pt);
	\node[Purple] at (2.2,1) {\small$A(S)$};
	\node[TealBlue] at (2.2,.6) {\small$P(S)$};
\end{tikzpicture}
\end{center}
\vfill
\newpage
More generally,
\begin{adjustwidth}{4mm}{4mm}
\textbf{Question B'.} Can we construct virtually $s$-step nilpotent groups with polynomial geodesic growth on top of any filiform nilpotent groups of type I
\[ \funcal{F}_s = \la y, z_1,z_2,\ldots,z_s \;\big|\; [y,z_i]=z_{i+1},\; [y,z_s]=[z_i,z_j]=1 \ra \;\text? \]
\end{adjustwidth}
Those groups have few abnormal curves, covering only the $z_1$-direction (see \cite{Nalon}).

\bigbreak

If Conjecture A holds, then Theorem \ref{thm:crit_step2} reduces the characterization of groups of polynomial geodesic growth (for symmetric generating set) to the following question:
\begin{adjustwidth}{4mm}{4mm}
\textbf{Problem C.} Characterize finite subgroups $F\le \GL_d(\ZZ)$ for which \vspace*{-1mm}
\begin{itemize}
	\item[\hyperlink{symm_cond}{(ii')}] there exists a finite \emph{symmetric} set $A\subset \ZZ^d$ such that $P=\Conv(A^F)$ is full-dimensional and each facet of $P$ contains at most one point of $A$.
\end{itemize}
\end{adjustwidth}
Ideally, the characterization should be algorithmic. Given $F$ as a finite set of integer-valued matrices, one should be able to decide whether or not the condition is satisfied. In contrast, \cite[Theorem 1]{bridson2012groups} quoted in the introduction gives the sufficient condition
\begin{itemize}[leftmargin=14.5mm]
	\item[(iii')] There exists $a\in\ZZ^d$ such that $P=\Conv(\{\pm a\}^F)$ is full-dimensional.
\end{itemize}
This condition is clearly algorithmic. It is therefore natural to ask
\begin{adjustwidth}{4mm}{4mm}
\textbf{Question D.} Are conditions (ii') and (iii') equivalent?
\end{adjustwidth}

\bigbreak

Finally, we would like to reiterate
\begin{adjustwidth}{4mm}{4mm}
	\textbf{Question E.} (\cite[Problem 13]{GrigorchukQ}) Does there exist a pair $(G,S)$ of intermediate geodesic growth and intermediate volume growth? 
\end{adjustwidth}
Exponential geodesic growth has been established in several examples of intermediate volume growth by \cite{bronnimann2016geodesic}, for instance the first Grigorchuk group with standard generating set. The simplest open example is the Fabrikowski-Gupta group. We isolate this question since an answer in either direction would require new insight on these groups.

\printbibliography
\end{document}